\documentclass[12pt,leqno]{amsart}
\usepackage{times}
\usepackage{amsmath,amssymb,amscd,amsthm}
\usepackage{latexsym}
\usepackage{epsfig}

\newtheorem{theorem}{Theorem}
\newtheorem{lemma}{Lemma}
\newtheorem{proposition}{Proposition}

\newtheorem{corollary}{Corollary}

\newtheorem{claim}{Claim}

\newcommand{\q}{\quad}
\newcommand{\qq}{\quad\quad}

\newcommand{\norm}[2]{{\left\| #1 \right\|}_{#2}}
\newcommand{\f}[2]{\frac{#1}{#2}}
\newcommand{\dpr}[2]{\langle #1,#2 \rangle}

\newcommand{\al}{\alpha}
\newcommand{\be}{\beta}
\newcommand{\ga}{\gamma}

\newcommand{\de}{\delta}

\newcommand{\De}{\Delta}
\newcommand{\ve}{\varepsilon}

\newcommand{\la}{\lambda}

\newcommand{\si}{\sigma}

\newcommand{\vp}{\varphi}
\newcommand{\om}{\omega}
\newcommand{\Om}{\Omega}

\newcommand{\rn}{{\mathbf R}^n}

\newcommand{\rone}{\mathbf R^1}
\newcommand{\rtwo}{\mathbf R^2}

\newcommand{\sn}{\mathbf S^{n-1}}
\newcommand{\sone}{\mathbf S^1}

\newcommand{\cf}{\mathcal F}

\newcommand{\cz}{\mathcal Z}

\newcommand{\cc}{\mathcal C}

\newcommand{\intl}{\int\limits}

\newcommand{\suml}{\sum\limits}
\newcommand{\supl}{\sup\limits}

\newcommand{\p}{\partial}

\newcommand{\beq}{\begin{equation}}
\newcommand{\eeq}{\end{equation}}
\newcommand{\beqna}{\begin{eqnarray*}}
\newcommand{\eeqna}{\end{eqnarray*}}
\newcommand{\beqn}{\begin{equation*}}
\newcommand{\eeqn}{\end{equation*}}
\newcommand{\bp}{\begin{proof}}
\newcommand{\ep}{\end{proof}}
\newcommand{\bprop}{\begin{proposition}}
\newcommand{\eprop}{\end{proposition}}
\newcommand{\bt}{\begin{theorem}}
\newcommand{\et}{\end{theorem}}
\newcommand{\bex}{\begin{Example}}
\newcommand{\eex}{\end{Example}}
\newcommand{\bc}{\begin{corollary}}
\newcommand{\ec}{\end{corollary}}
\newcommand{\bcl}{\begin{claim}}
\newcommand{\ecl}{\end{claim}}
\newcommand{\bl}{\begin{lemma}}
\newcommand{\el}{\end{lemma}}

\begin{document}

\title
[Pseudodifferential operators with rough symbols]
{Pseudodifferential operators with rough symbols}

\author{Atanas Stefanov}

\address{
Department of Mathematics,
University of Kansas,
Lawrence, KS~66045, USA}
\date{\today}

\thanks{Supported in part by  NSF-DMS 0300511}



\begin{abstract}
In this work, we develop  $L^p$ 
boundedness theory for pseudodifferential operators 
with rough (not even continuous in general) symbols 
in the $x$ variable. Moreover, the  
$B(L^p)$ operator norms are estimated explicitly in 
terms of scale invariant quantities involving the symbols.  All 
the estimates are shown to be sharp with respect to the required 
smoothness in the 
$\xi$ variable. As a corollary, we obtain $L^p$  bounds for (smoothed out
versions of) the maximal directional Hilbert transform 
and the Carleson operator. 
\end{abstract}

\maketitle
\date{today}

\section{Introduction}
In this paper, we are concerned with the $L^p$ mapping properties of 
the pseudodifferential operators in the form 
\begin{equation}
\label{eq:1}
T_\si f(x)=\intl_{\rn} \si(x,\xi) e^{2\pi i \xi x} \hat{f}(\xi) d\xi.
\end{equation} 
The operators $T_\si$  have been subject of continuous interest 
since the sixties. We should mention that their usefullness in 
 the study of partial differential equations have been realized  much 
earlier, but it seems that their systematic study began with the 
fundamental works of Kohn and Nirenberg, \cite{KN} and  H\"ormander, 
\cite{Hor}. 

To describe the 
results obtained in these early papers, 
define the H\"ormander's class $S^m$, which consists 
of all functions $\si(x,\xi)$, so that 
\begin{equation}
\label{eq:2}
|D^\beta_x D^\al_\xi \si(x, \xi)|\leq C_{\al, \be} (1+|\xi|)^{m-|\al|}.
\end{equation}
for all multiindices $\al, \be$. 
A classical theorem in \cite{Hor} then states that 
$Op(\si):H^{s+m,p}\to H^{s, p}$  for all $s\geq 0$ and 
$1<p<\infty$. In particular, $Op(\si):L^p\to L^p$,  $1<p<\infty$, whenever the symbol $\si\in S^{m}$. 
Subsequent improvements of these methods established the 
boundedness of $Op(\si)$ (basically under the assumption $\si\in S^m$ for 
appropriate $m$) 
 to various related function spaces, like 
Besov, Triebel-Lizorkin spaces to name a few, but we will not review 
those here, since they fall  
outside of the scope of this paper. 

It is worth mentioning however, that the simple to verify condition 
\eqref{eq:2} is the one arising in many  applications. 
The $L^2$ boundedness plays special role in the theory and that is why
we discuss it separately. 

The class of symbols $S^m_{\rho, \delta}$, 
defined via 
\begin{equation}
\label{eq:3}
|D^\beta_x D^\al_\xi \si(x, \xi)|\leq C_{\al, \be} 
(1+|\xi|)^{m-\rho|\al|+\delta|\beta|}.
\end{equation}
represents a larger set of symbols than $S^m=S^m_{1, 0}$, which has 
subsequently found 
 applications in  local 
solvability for linear PDE's, \cite{BF}.

Here, we have to mention the celbrated  result of 
Calder\'on-Vaillancourt, 
\cite{CV1}, \cite{CV2} which  states  that $L^2$ boundedness for $T_\si$ 
holds, whenever  
$\si\in S^{0}_{\rho, \rho}$, $0\leq \rho<1$, 
 whereas $S^{0}_{1,1}$ is a ``forbidden'' class, 
in the sense that there are symbols in that class, 
which give rise to unbounded on 
$L^2$ operators. We should mention here the work of Cordes \cite{Cor}, 
 who improved the 
result for $S^0_{0, 0}$ by requiring that \eqref{eq:3} holds only for 
$|\al|, |\beta|\leq [n/2]+1$. 

Regarding less regular in $x$ symbols,  
for any modulus of continuity 
$\om:R_+\to R_+$ (that is,  an  increasing and  continuous  function), 
define the space $C^\om$ of all uniformly continuous and bounded functions 
$u:\rn\to \cc$,  satisfying 
$$
|u(x+y)-u(x)|\leq \om(|y|). 
$$
The following class of symbols was introduced and 
 studied by  Coifman-Meyer, \cite{CM}. More precisely, 
let $\si(x,\xi)\in C^\om S^{0}_{1,0}$, which means  it  satisfies 
$$\sup_{x, y, \xi} <|\xi|>^{|\al|} 
|D_\xi^\al [\si(x+y, \xi)-\si(x, \xi)]| \leq C_\al \om(|y|)
$$
and assume that $\sum_{j>0} \om(2^{-j})^2<\infty$. Then for 
all $1<p<\infty$, $
Op(\si):L^p\to L^p. 
$
The condition $\sum_j \om(2^{-j})^2<\infty$   is clearly very mild continuity assumption for the function $x\to \si(x, \xi)$. 
In particular, one sees that $\cup_{\ga>0} 
C^\ga S^0\subset C^\om S^{0}_{1,0}$. 
Related results can be found in the work of M. Taylor, \cite{Taylor} (see Proposition 2.4, p. 23) and J. Marschall, \cite{Marschall} where the 
spaces $C^\om$ are replaced by $H^{\ve,p}$ spaces with $p$ as large as one wish and $0<\ve=\ve(p)<<1$  (see also \cite{Taylor}, p. 61)

One of the purposes of this work is  to get away from the 
continuity requirements  on $x\to \si(x, \xi)$. Even more importantly, 
we would like  to replace the pointwise conditions on the derivatives of $\xi$ 
by averaged ones. This particular point has not been thoroughly 
explored appropriately  
 in the literature in the author's opinion, 
see Theorem \ref{theo:1} below.  

On the other hand, a particular motivation for such considerations is 
provided by the recent  papers of Rodnianski-Tao \cite{RT} and 
the author  \cite{Stefanov}, where concrete parametrices 
(i.e. pseudodifferential operators,  representing 
approximate solutions to certain PDE's) 
were constructed for the 
solutions of certain first order perturbation of the wave 
and Schr\"odinger equations.  
A very quick inspection of these examples shows that\footnote{Most 
readers are 
 likely to have their own  fairly long list 
with  favorite examples, for which 
 the H\"ormander condition fails.} 
 {\it 
they do not obey
pointwise conditions on the derivatives on the 
symbols}   and thus, these methods fail 
to imply  $L^2$ bounds for 
these (and related problems). Moreover, one often times has to deal with the 
situation, where the maps $\xi\to \si(x, \xi)$ 
are not smooth in a pointwise sense. 
On the other hand, one may still be able to control averaged quantities  like 
\begin{equation}
\label{eq:5}
\sup_x \norm{\si(x, \xi)}{H^{n/2}_\xi}<\infty.
\end{equation}
This will be our treshold condition for $L^2$ boundedness, 
which we try to achieve. \\ Heuristically at least, \eqref{eq:5} 
must be ``enough'' in some sense, 
 since if we had simple symbols like $\si(x,\xi)=\si_1(x) \si_2(\xi)$, 
then the  $L^2$ boundedness of $Op(\si)$ is equivalent to 
   $\norm{\si_1}{L^\infty_x}<\infty, 
\norm{\si_2}{L^\infty_\xi}<\infty$. Clearly, 
$\norm{\si_2}{L^\infty_\xi(\rn)}$ just fails to be controlled by 
\eqref{eq:5}, but on the other hand, the quanitity in \eqref{eq:5} 
is controlled by the 
appropriate Besov space  $B^{n/2}_{2,1}$ norm. 

 A final motivation for 
the current study is  to achieve a scale
 invariant condition, which gives an estimate of the 
$L^2\to L^2$ ($L^p\to L^p$) norm of 
$Op(\si)$ in terms of a {\it scale invariant quantity}, 
that is,  we aim at showing 
an estimate, 
$$
\norm{Op(\si)}{L^p\to L^p}\leq C \norm{\si}{Y} \norm{f}{L^p}, 
$$
where for every $\la\neq 0$, one has $\norm{\si(\la \cdot,\la^{-1} 
\cdot) }{Y}=\norm{\si}{Y}$. 

In that regard, note that the  condition (which is one of the requirements of 
 the  H\"ormander class $S^0$)  
\begin{equation}
\label{eq:6}
\sup_{x}| D_\xi^\al \si(x, \xi)|\leq  C_\al |\xi|^{-|\al|}
\end{equation}
is   scale invariant in the sense described above. Moreover, 
by the standard Calder\'on-Zygmund theory (see \cite{Stein}), 
the pointwise condition \eqref{eq:6} together with  
$\|T_\si\|_{L^2\to L^2}<\infty$  implies 
$$
T_\si f(x)=\int K(x, x-y) f(y) dy,
$$
where $K(x,\cdot)$ satisfies the H\"ormander-Mihlin conditions, namely $|K(x, z)|\leq C|z|^{-n}$ and 
$|\nabla_z K(x,z)|\leq C|z|^{-n-1}$, where the constant $C$ depends on 
the constants 
$C_\al: |\al|<[n/2]+1$ in \eqref{eq:6}.  
This  in turn is enough to conclude  
that $T_\si:L^p \to L^p$ for all $1<p\leq 2$ and in fact there is the 
endpoint estimate  $T_\si:L^1\to L^{1, \infty}$. 

\subsection{$L^p$ estimates for PDO with rough symbols - statement of results} 
We start now with our main theorems, which concern the  $L^2$ and the $L^p$ boundedness for pseudodifferential 
operators $Op(\si)$ with rough symbols. Our first result establishes that a Besov space version of \eqref{eq:5} is enough for $L^2$ boundedness and 
the result is sharp. 
\begin{theorem}($L^2$ bounds) 
\label{theo:1}
Let $\si(x, \xi):\rn\times \rn\to \cc$ and $T_\si$ is the 
corresponding pseudodifferential operator. Then 
\begin{equation}
\label{eq:7}
\norm{T_\si }{L^2(\rn)\to L^2(\rn)}\leq C (\suml_l 2^{l n/2} \supl_x \|P_l^\xi \si(x, \cdot)\|_{L^2(\rn)}), 
\end{equation} 
where $P_l^\xi$ is the Littlewood-Paley operator in the $\xi$ variable. 

Moreover, the result is sharp in the  following sense:  
for every $p>2$, 
there exists $\si(x,\xi)$ so that $\sup_{x} | D_\xi^{\al}\si(x, \xi)|\leq C_\al |\xi|^{-|\al|}$ and $\sup_x \norm{\si(x, \cdot)}{W^{p, n/p}}<\infty$,  
but  $T_\si$ fails to be bounded on $L^2(\rn)$. 
\end{theorem}
{\bf Remark:} 
\begin{enumerate} 
\item Note that the estimate on $T_\si$ is scale invariant. 
\item The sharpness claim of the theorem, roughly speaking, 
shows that in the scale of spaces\footnote{Note that 
these  spaces scale the same and moreover 
by Sobolev embedding these are strictly decreasing sequence, at least  
for $2\leq p<\infty$.} $W^{p, n/p}$, $\infty\geq p\geq 2$, one may 
not require anything less than $W^{2,n/2}=H^{n/2}$ of the symbol in 
order  to ensure $L^2$ boundedness.
\item The counterexample to which we refer in Theorem \ref{theo:1} 
is a simple variation of the well-known example of $\si\in 
S^{0}_{1,1}$, the ``forbidden class'', 
 which fails to be $L^2$ bounded, see \cite{Stein}, p. 272 and Section 
\ref{sec:counter} below. 
 \end{enumerate}
Our next result concerns $L^p$ boundedness for $T_\si$. 
\begin{theorem}($L^p$ bounds) 
\label{theo:3}
For the pseudodifferential operator $T_\si$ there 
is the estimate for all $2<p\leq \infty$, 
\begin{equation}
\label{eq:22}
\norm{T_\si }{L^p(\rn)\to L^p(\rn)}\leq 
C (\suml_l 2^{l n/2} \supl_x \|P_l^\xi \si(x, \cdot)\|_{L^2(\rn)}), 
\end{equation}
For the range $1<p<2$ and indeed for the weak type $(1,1)$, there is 
\begin{equation}
\label{eq:8}
\norm{T_\si}{L^p\to L^p}+\norm{T_\si}{L^1\to L^{1, \infty}}\leq C 
\suml_l  2^{ l n} \supl_{x}  \|P_l^\xi \si(x, \cdot)\|_{L^1(\rn)}. 
\end{equation}
Alternatively, if one assumes the $L^2$ bound, together with 
\eqref{eq:6}, one still gets $L^p\to L^p$, $1<p\leq 2$, and 
in fact weak type $(1,1)$ 
  bounds. Moreover, 
$$
\norm{T_\si}{L^p\to L^p}\leq C (
\suml_l  2^{ l n/2} \supl_{x}  \|P_l^\xi \si(x, \cdot)\|_{L^2(\rn)}+
\sup_{|\al|<[n/2]+1} \sup_{x,\xi} |\xi|^{|\al|} | D_\xi^\al \si(x, \xi)|).
$$
\end{theorem}
As we pointed out in Theorem \ref{theo:1}, the estimates are 
essentially sharp for $L^p$,  $2\leq p<\infty$ boundedness. 
The following  corollary gives even more precise 
condition under which a symbol $\si$ will give rise to a bounded 
operator on $L^q$ in the case of a given $1<q<2$
\begin{corollary}
\label{cor:1}
Let $1<q<2$.  Then 
$$
\norm{T_\si}{L^q\to L^q}\leq C\suml_l 2^{ln/q} \supl_x 
\|P_l^\xi \si(x, \cdot)\|_{L^{q}(\rn)}. 
$$
\end{corollary}
Clearly the proof follows by  interpolation 
from  the $L^2$ estimates  in Theorem \ref{theo:1} and the weak type $(1,1)$ 
estimates of Theorem \ref{theo:3}.

\subsection{PDO's with  homogeneous of degree 
zero symbols - statement of  results}

Regarding symbols that are  homogeneous of degree zero, 
i.e. $\si(x, \xi)=q(x, \xi/|\xi)$, where $q:\rn\times 
\sn\to \cc$, we  obtain more precise results in terms of 
the smoothness of $q$. \\
Note that the classical H\"ormander condition  requires 
pointwise smoothness of the function $q$ in both variables. 
Our result on the other hand requires much less than that.  
\begin{theorem}($L^p$ bounds for  homogeneous of degree zero symbols)\\
\label{theo:4}
Let $q:\rn\times \sn\to \cc$. Let  
$$
T_q f(x)= \intl_{\rn} q (x,\xi/|\xi|) e^{2\pi i \xi x} 
\hat{f}(\xi) d\xi.
$$
Then $T_q :L^2\to L^2$, if $\sum_l 2^{l(n-1)/2} 
\|P_l^{\xi/|\xi} q\|_{L^2(\sn)}<\infty$ and in fact 
\begin{equation}
\label{eq:9}
\norm{T_q}{L^2\to L^2}\leq C \sum_l 2^{l(n-1)/2} 
\sup_x \|P_l^{\xi/|\xi|} q(x, \cdot)\|_{L^2(\sn)}.
\end{equation}
Concerning $L^p$ bounds, we have for every  $2\leq p\leq  \infty$. 
\begin{equation}
\label{eq:10}
\norm{T_q}{B^0_{p,1}\to L^p}\leq C_n (\sum_l 2^{l(n-1)/p'} 
\sup_x \|P_l^{\xi/|\xi|} q\|_{L^{p'}(\sn)}).
\end{equation}
Note that in \eqref{eq:10}, 
the constant $C_n$  is independent of  $p,r$. 
\end{theorem} 
{\bf Remark:}
\begin{enumerate}
\item It would be interesting to see whether the usual $L^p\to L^p$  
 boundedness  holds true. 
\item 
Note that there is no  weak type $(1,1)$ statement in 
 Theorem \ref{theo:4}. This is a  difficult 
issue even for multipliers. 
\end{enumerate}
\noindent The sharpness statement associated with Theorem \ref{theo:4} is 
\begin{proposition}
\label{prop:5}
For every $N>1$, there exists a homogeneous of degree zero  symbol 
$\si(x, \xi):\rtwo\times \rtwo\to \rone$, 
so that $\sup_{x, \xi} |\si(x, \xi)|<\infty$ and 
$\sup_x \norm{\si(x,\xi)}{W^{1, 1}(\sone)}<\infty$, 
and so that $\norm{T_\si}{L^2\to L^2}>N$. 
\end{proposition}
The counterexample considered here is a smoothed out version  of the 
maximal directional Hilbert transform in the plane $H_* f(x)=\sup_{u\in \sone} 
|H_u f(x)|$. We mention the spectacular recent result of 
Lacey and Li, \cite{LL} showing the boundedness of  
$H_*$ on $L^p(\rtwo): 2<p<\infty$ with a 
$H_*: L^2(\rtwo)\to L^{2, \infty}(\rtwo)$ as an endpoint 
estimate. Note that the $L^2\to L^2$ bound fails, as elementary examples show, see \cite{LL}.  
We verify later that the condition \eqref{eq:10} 
just fails  for the (smoothed out) multiplier $\si$  of $H_*$  in two 
dimensions, 
but on the other hand the condition 
$\norm{\si(x,\xi)}{W^{1, 1}(\sone)}<\infty$ holds. 

This example will show that the Besov spaces requirements for $\si$ in 
\eqref{eq:9} and \eqref{eq:10} cannot be replaced by 
Sobolev spaces and/or spaces with less derivatives.  

\subsection{PDO's with  radial symbols.}
Finally, we consider the case of radial symbols. That is for 
$\rho:\rn\times R^1_+\to \cc$ and 
$$
T_\rho f(x)=\int_{\rn} \rho(x, |\xi|) e^{2\pi i \xi x} \hat{f}(\xi) d\xi.
$$
\begin{theorem}
\label{theo:radial}
The operator $T_\rho:L^2\to L^2$, if $\sum_l 2^{l/2} \supl_x \|P^{|\xi|}_l
\rho(x, \cdot)\|_{L^2(\rone)}<\infty$. In fact, 
$$
\|T_\rho\|_{L^2\to L^2}\leq C \sum_l 2^{l/2} \supl_x \|P^{|\xi|}_l
\rho(x, \cdot)\|_{L^2(\rone)}.
$$

\end{theorem}
Clearly, establishing  $L^p, p\neq 2$ bounds for simple radial {\it symbols} 
is already a notoriously difficult problem. 
One only needs to point out to the Bochner-Riez multiplier 
$(1-|\xi|^2)^{\de}_+$ (which satisfy $L^p$ bounds  only in certain range of
$p's$, 
depending on the dimension and
$\de$) or even the simpler ``thin annulus'' multiplier 
$\vp(2^m(1-|\xi|^2))$ to understand the difficulty of the problem in general. 

\section{Applications}
In this section, we demonstrate the effectiveness of the $L^p$ boundedness
theorems for rough PDO's. We will mostly concentrate on application to 
maximal functions and operators\footnote{In addition, the author has also 
identified   several applications to bilinear/multilinear  operators of 
  importance to certain
dispersive PDE's, which  will be addressed in a future publication.} 
Some of our examples will be well-known results for  maximal operators, 
while others will be a higher dimensional extensions of such results. 

\subsection{Almost everywhere convergence for Cesaro sums of $L^p$ functions  in 1 D} 
We start with Cesaro's sum for Fourier series in one space dimension. For any
$\de>0$, define 
$$
\cc_\de f(x)=\sup_{u>0} \int_{\rone} (1-\xi^2/u^2)^{\de}_+  
e^{2\pi i \xi x} \hat{f}(\xi) d\xi.
$$ 
Clearly, as a limit as $\de\to 0$, we get the Carleson's operator. 
Unfortunately, one cannot  conclude that 
$\sup_{\de>0}\|C_\de\|_{L^p}<\infty$, for
that would imply the famous Carleson-Hunt theorem. 
On the other hand, define the maximal ''thin
interval operator''
$$
T_m f(x)=\sup_{u>0} \int_{\rone} \vp(2^m(1-\xi^2/u^2))  
e^{2\pi i \xi x} \hat{f}(\xi) d\xi.
$$
A 
simple argument based on (the proof of) Theorem \ref{theo:3} yields 
\begin{proposition}
\label{prop:Carl}
For any $\ve>0, 1<p<\infty$, there exists $C_{p, \ve}$, so that 
\begin{equation}
\label{eq:908}
\sup_m \|T_m \|_{L^p(\rone)\to L^p(\rone)}\leq C_p. 
\end{equation}
In fact, there is the more general pointwise 
bound $\supl_m |T_m f(x)|\leq C M (\sup_k |P_k f|)(x)$, 
which implies \eqref{eq:908} as well as 
\begin{equation}
\label{eq:909}
\|C_\de\|_{L^p\to L^p}+
\|C_\de\|_{F^0_{1, \infty}(\rone)\to L^{1,\infty}(\rone)}\leq C_{\de,p}. 
\end{equation}
\end{proposition}
{\bf Remark:} 
\begin{itemize}
\item Note that this result, while clearly inferior to the 
Carleson-Hunt 
theorem still implies a.e. convergence for any Cesaro summability method, when
applied to $L^p(\rone)$ functions, and in fact for the larger class of $F^0_{p, \infty}$ functions.  
\item Using the method of proof here, one may actually prove $L^2$ estimates\footnote{And in fact
$L^p\to L^p$ estimates for $p=p(\de)$ close to $2$.} for the {\it maximal
Bochner-Riesz operator} 
$$
BR_\de f(x)=\sup_{u>0} \int_{\rn} (1-|\xi|^2/u^2)^{\de}_+  
e^{2\pi i \xi x} \hat{f}(\xi) d\xi.
$$
{\it in any dimension}. 
\end{itemize}
\begin{proof}
It clearly suffices to show the pointiwise estimate 
$|T_m f(x)|\leq C M (\sup_k P_k f)(x)$ for any $k$.
The  statements about $B^0_{p,1}\to L^p$ bounds  follow by elementary Littlewood-Paley theory and
the $l^p$ bounds for the Hardy-Littlewood maximal function. The restricted-to-weak estimate 
$F^0_{1,\infty}\to L^{1,\infty}$ for $C_\de$ follows by summing an 
{\it exponentially decaying} series in
the quasi-Banach space $L^{1, \infty}$. 

By support considerations, it is clear that 
\begin{eqnarray*}
T_m f(x) & = &\suml_k \sup_{u>0} \int_{\rone} \vp(2^m(1-\xi^2/u^2))  
e^{2\pi i \xi x}\vp(2^{-k}\xi) \hat{f}(\xi) d\xi= \\
&=&
\suml_k \sup_{u\in (2^{k-2}, 2^{k+2})} \int_{\rone} \vp(2^m(1-\xi^2/u^2))  
e^{2\pi i \xi x}\vp(2^{-k}\xi) \hat{f}(\xi) d\xi= \\
&=& \suml_k T_{m, u(\cdot)\in (2^{k-2}, 2^{k+2})} f_k.
\end{eqnarray*}
Clearly, the requirement $u\in (2^{k-2}, 2^{k+2})$ creates (almost) 
disjointness in the $x$ support, whence 
\begin{equation}
\label{eq:819}
|T_m f(x)|\leq C \supl_k |T_{m, u(\cdot)\in (2^{k-2}, 2^{k+2})} f_k(x)|.
\end{equation}
Our basic claim is that 
\begin{equation}
\label{eq:820}
|T_{m, u(\cdot)\in (2^{k-2}, 2^{k+2})} f_k(x)|\leq C M(f_k).
\end{equation}
Clearly \eqref{eq:819} and \eqref{eq:820} imply  $\sup_m |T_m f(x)|\leq C M(\supl_k |f_k|)$, 
whence the Proposition \ref{prop:Carl}.  \\
By scale invariance, \eqref{eq:820} reduces to the case $k=0$, that is we need to show 
$$
|T_{m, u(x)\in (1/4, 4)} P_0 f(x)|\leq C M (P_0 f)(x).
$$
for any Schwartz function $f$ and any 
$m>>1$. By \eqref{eq:32} (in the proof of Theorem \ref{theo:3} below), it will suffice to
show 
\begin{equation}
\label{eq:920}
\suml_l 2^l \supl_x \|P_l^\xi [\vp(2^m(1-\xi^2/u(x)^2)) \vp(\xi)]\|_{L^1(\rone)}\lesssim 1.
\end{equation}
for any measurable function $u$, which takes its values in $(1/4, 4)$.  \\
For \eqref{eq:920}, we have 
\begin{eqnarray*}
& & \suml_{l<m}  2^l \supl_x \|P_l^\xi [\vp(2^m(1-\xi^2/u(x)^2)) \vp(\xi)]\|_{L^1_\xi(\rone)} \\ 
& & \leq \suml_{l<m} 2^l \supl_{u\in (1/4, 4)} \int|\vp(2^m(1-\xi^2/u^2)) \vp(\xi)| d\xi\lesssim 
\suml_{l<m} 2^{l-m}\lesssim 1, 
\end{eqnarray*}
while 
\begin{eqnarray*}
& & \suml_{l\geq m}  2^l \supl_x 
\|P_l^\xi [\vp(2^m(1-\xi^2/u(x)^2)) \vp(\xi)]\|_{L^1_\xi(\rone)} \\
& &\leq \suml_{l\geq m}  2^{-l} 
\supl_{u\in (1/4, 4)} \int|\f{d^2}{d \xi^2} \vp(2^m(1-\xi^2/u^2)) \vp(\xi)|d\xi\leq C 
\suml_{l\geq m} 2^{-l+m}\lesssim 1.
\end{eqnarray*}
\end{proof}
\subsection{Maximal directional Hilbert transforms and the Kakeya  maximal
function.}
Another interesting application is provided by the 
directional Hilbert transform
in dimensions  $n\geq 2$. Namely, take 
$$
H^*_\de=\sup_{u\in \sn} \int (\dpr{u}{\xi/|\xi|})^{\de}_+ 
e^{2\pi i \xi x} \hat{f}(\xi)\vp(\xi) d\xi,
$$
where $supp \vp\subset \{1/2<|\xi|<2\}$. \\
As $\de\to 0$, we obtain the operator $f\to 
\sup_{u\in \sn} \int_{\{\dpr{u}{\xi}>0\}}
e^{2\pi i \xi x} \hat{f}(\xi) d\xi$, which is closely related to the maximal
directional Hilbert transform  
$$
H_* f(x)=\sup_{u}|H_u f(x)|=\sup_u |\int sgn(u, \xi) 
e^{2\pi i \xi x} \hat{f}(\xi) \vp(\xi)d\xi|.
$$
$H_*$ was of course shown to be $L^p(\rtwo), p>2$ bounded by 
Lacey and Li, \cite{LL} by very sophisticated time-frequency analysis methods. 
\begin{proposition}
\label{prop:mdh}
For the ``thin big circle'' multiplier 
$$
T_m f(x)= \sup_{u\in \sn} |\int_{\rn} \vp(2^m\dpr{u}{\xi/|\xi|}) 
e^{2\pi i \xi x} \hat{f}(\xi) \vp(\xi)d\xi|.
$$
we have 
\begin{equation}
\label{eq:kak}
\|T_m f\|_{L^2\to L^2}\leq C_\ve  2^{m(n/2-1)}
\end{equation}
In particular 
$$
\|H^*_\de\|_{L^2(\rn)\to L^2(\rn)}\leq C_{p, \ve, \de} 2^{n/2-1}.
$$
\end{proposition}
{\bf Remark:} 
\begin{itemize}
\item We believe that  the operator $T_m$ ($m>>1$) 
has a particular connection to the 
Kakeya maximal function and the corresponding Kakeya problem. 
Indeed, the kernel
of the corresponding singular integral behaves like  a 
($L^1$ normalized) characteristic function of a rectangle with long side along $u$ of length $2^m$
and $(n-1)$ short sides of length $1$ in  the transverse directions! 
\item In relation to that, one expects the conjectured Kakeya bounds 
$$
\|T_m f\|_{L^p\to L^p}\leq C_\ve   2^{m(n/p-1)}
$$
for $p\leq n$ 
to hold, while one only gets 
$$
\|T_m f\|_{L^p\to L^p}\leq C_\ve   2^{m(n/p-2/p)}
$$
as a consequence of Theorem \ref{theo:4}. Nevertheless, the two match when 
$p=2$. So it seems that  \eqref{eq:kak}, at
least in principle,  captures the Kakeya conjecture for 
$p=2$ in general and in
particular the full Kakeya conjecture in two dimensions. 

Since our estimates do not seem to contribute much toward the resolution of any
new Kakeya estimates, we do not pursue here the exact relationship between $T_m$
and the Kakeya maximal operator, although from our heuristic arguments above it
should be clear that it is a close one. 
\end{itemize}
\begin{proof}
We proceed as in the proof of Proposition \ref{prop:Carl}. We need only show 
\begin{equation}
\label{eq:990}
\suml_{l} 2^{l(n-1)/2} \supl_x \|P_l^{\xi/|\xi|} 
\vp(2^m\dpr{u(x)}{\xi/|\xi|})\|_{L^2(\sn)}\lesssim 1.
\end{equation}
We have 
\begin{eqnarray*}
& & \suml_{l<m} 2^{l(n-1)/2} \supl_x \|P_l^{\xi/|\xi|} 
\vp(2^m\dpr{u(x)}{\xi/|\xi|})\|_{L^2(\sn)}\\ 
& &  \leq C \suml_{l<m} 2^{l(n-1)/2}\sup_u 
\| \vp(2^m\dpr{u(x)}{\xi/|\xi|})\|_{L^2(\sn)}\\
& &\leq  C \suml_{l<m} 2^{l(n-1)/2} 2^{-m/2}
\lesssim 2^{m(n/2-1)}.
\end{eqnarray*}
\begin{eqnarray*}
& & \suml_{l\geq m} 2^{l(n-1)/2} \supl_x \|P_l^{\xi/|\xi|} 
\vp(2^m\dpr{u(x)}{\xi/|\xi|})\|_{L^2(\sn)}\\ 
& & C \leq 
\suml_{l\geq m} 2^{-l(n-1)/2} \supl_u \|\Om^{(n-1)}
\vp(2^m\dpr{u(x)}{\xi/|\xi|})\|_{L^2(\sn)} \\
& & \leq C  \suml_{l\geq m} 2^{-l(n-1)/2} 2^{m(n-1)} 2^{-m/2}\lesssim 2^{m(n/2-1)}. 
\end{eqnarray*}
\end{proof}

\subsection{Estimates on $T_{\si_1 \si_2}$, $T_{e^{ \si}}$ etc.}
We now present a result, which allows us to treat pseudodifferential 
 operators, whose symbols are products, exponentials (or more generally entire functions) of symbols, which satisfy the requirements in 
Theorems \ref{theo:1}, \ref{theo:3}, \ref{theo:4}. We would like to point out that similar in 
spirit (by essentially requiring $n/2+\ve$ derivatives in $L^2$, but in a more general setting) 
functional calculus type result was obtained in \cite{SeegerSogge}. 
\begin{proposition}
\label{prop:7}
Let $\si, 
\si_1, \si_2:\rn\times \rn\to \cc$, so that $T_{\si_j}, j=1,2$ define 
$L^p$ bounded operators, as in \eqref{eq:7}, \eqref{eq:8}.  
Then for every $1<p<\infty$, $T_{\si_1\si_2}$ is 
 $L^p$ bounded and 
\begin{eqnarray}
\label{eq:814}
& & \norm{T_{\si_1 \si_2}}{L^2(\rn)\to L^2(\rn)}\leq C  \prod_{j=1}^2 
(\suml_l 2^{l n/2} \supl_x \|P_l^\xi \si_j(x, \cdot)\|_{L^2(\rn)}) , \\
\label{eq:815}
& & \norm{T_{\si_1 \si_2}}{L^p(\rn)\to L^p(\rn)}\leq C 
\prod_{j=1}^2 
(\suml_l 2^{l n} \supl_x \|P_l^\xi \si_j(x, \cdot)\|_{L^1(\rn)}).
\end{eqnarray}
In the same spirit, $T_{e^{\si}}$ is also $L^p$ bounded and there is 
\begin{equation}
\label{eq:816}
\norm{T_{e^\si}}{L^p(\rn)\to L^p(\rn)}\leq C  
exp(\suml_l 2^{l n} \supl_x \|P_l^\xi \si(x, \cdot)\|_{L^1(\rn)}) 
\end{equation}
Similar statements can be made for homogeneous of 
degree zero symbols \\ $\mu_1(x, \xi/|\xi|), \mu_1(x, \xi/|\xi|)$,  
\begin{eqnarray}
\label{eq:817}
& & 
\norm{T_{\mu_1\mu_2}}{L^2(\rn)\to L^2(\rn)}\leq C  
\prod_{j=1}^2 (\suml_l 2^{l (n-1)/2} \supl_x \|P_l^{\xi/|\xi|} \mu_j(x, 
\cdot)\|_{L^1(\rn)}), \\
& & 
\label{eq:818}
\norm{T_{e^\mu}}{L^2(\rn)\to L^2(\rn)}\leq C  
exp(\suml_l 2^{l (n-1)} \supl_x \|P_l^{\xi/|\xi|} \mu(x, 
\cdot)\|_{L^1(\rn)}).
\end{eqnarray}

\end{proposition}
The proof of Proposition \ref{prop:7} is based on the corresponding 
Theorem for $L^p$ boundedness, 
combined with the fact that our requirements form a Banach 
algebra under the multiplication. Take for example \eqref{eq:814}. 
By Theorem \ref{theo:1}, we have 
$$
\norm{T_{\si_1 \si_2}}{L^2(\rn)\to L^2(\rn)}\leq C   
(\suml_l 2^{l n/2} \supl_x \|P_l^\xi (\si_1\si_2)(x, \cdot)\|_{L^2(\rn)})
$$
We finish by invoking the estimate
\begin{equation}
\label{eq:900}
\suml_l 2^{l n/2} \supl_x \|P_l^\xi (\si_1\si_2)(x, \cdot)\|_{L^2(\rn)}\leq 
C \prod_{j=1}^2 
(\suml_l 2^{l n/2} \supl_x \|P_l^\xi \si_j(x, \cdot)\|_{L^2(\rn)}).
\end{equation}
where this last inequality essentially means that $B^{n/2}_{2,1}$ 
is a Banach algebra of functions\footnote{This is well-known, but 
can be verified easily by 
means of the Kato-Ponce estimate 
$\|\p^{n/2} (u v)\|_{L^2}\leq C (\|\p^{n/2} u\|_{L^2}\norm{v}{L^\infty}+ 
\|\p^{n/2} v\|_{L^2}\norm{u}{L^\infty})$, the embedding 
$B^{n/2}_{2,1}\hookrightarrow L^\infty$ 
and  some Littlewood-Paley theory.}.

The argument above can be performed for the proof  of 
\eqref{eq:815}. For  \eqref{eq:816} (and more generally for 
any symbols of the form $g(\si)$, where $g$ is entire function), 
one  iterates  the product estimate 
\eqref{eq:900} to 
\begin{equation}
\label{eq:901}
\suml_l 2^{l n/2} \supl_x \|P_l^\xi (e^\si)(x, \cdot)\|_{L^2(\rn)}\leq 
C exp(\suml_l 2^{l n/2} \supl_x \|P_l^\xi \si(x, \cdot)\|_{L^2(\rn)}).
\end{equation}
For the proof of \eqref{eq:817}, \eqref{eq:818}, one has to use 
 the fact that  $B^{(n-1)/2}_{2, 1}(\sn)$ is a 
Banach algebra as well, whence one gets an estimate similar to  
\eqref{eq:900} and \eqref{eq:901}.

\section{Preliminaries}
\label{sec:prelim}
\noindent We start by introducing some basic concepts in Fourier analysis. 
\subsection{Fourier analysis on $\rn$}
First, define the Fourier transform and its inverse  
\begin{eqnarray*}
& &
\hat{f}(\xi)=\intl_{\rn} f(x)e^{-2\pi i x\cdot  \xi} dx, \\
& &
 f(x)=\intl_{\rn} \hat{f}(\xi)e^{2\pi i x\cdot  \xi} d\xi. 
\end{eqnarray*}
For a 
positive, smooth and  even function $\chi:\rn\to \rone_+$, supported in 
 $\{\xi:|\xi|\leq 2\}$ and so that $\chi(\xi)=1$ for all
$|\xi|\leq 1$.  
 Define $\vp(\xi)=\chi(\xi)-\chi(2\xi)$, 
which is supported in the annulus $1/2\leq |\xi|\leq 2$. Clearly 
$\sum_{k\in \cz} \vp(2^{-k} \xi)=1$ for 
all\footnote{The discussion henceforth will be 
for $\rn$, unless explicitely specified otherwise.}  
$\xi\neq 0$. 

The $k^{th}$ Littlewood-Paley projection is given by 
$\widehat{P_k f}(\xi)=\vp(2^{-k}\xi) \hat{f}(\xi)$.  
Note that the kernel of $P_k$ is  integrable uniformly  in $k$  and thus 
$P_k :L^p\to L^p$ for $1\leq p\leq \infty$ and  
$\norm{P_k}{L^p\to L^p}\leq 
C_n \norm{\hat{\chi}}{L^1}$. 
In particular, the  bounds are independent of $k$. 

It is a standard observation that $\nabla P_k = P_k \nabla = 
2^{k}\tilde{P_k}$, where $\tilde{P_k}$ is a multiplier type 
operator similar to $P_k$ and thus 
$\norm{P_k \nabla  \psi}{L^p}\sim 2^k \norm{P_k \psi}{L^p}$.  
In what follows, we will use the notation 
$P_l^x f(x, \xi)$ to denote Littlewood-Paley operator acting on 
the variable $x$, and $P_l^\xi f(x, \xi)$ will be a Littlewood-Paley operator in the variable $\xi$. That is 
\begin{eqnarray*}
& & 
P_l^x f(x, \xi)=2^{l n} \int \hat{\vp}(2^l(x-y)) f(y, \xi) dy, \\
& & P_l^\xi f(x, \xi)=2^{l n} 
\int \hat{\vp}(2^l(\xi-\eta)) f(x, \eta) d\eta
\end{eqnarray*}
The  Bernstein inequality takes the form
$$
\norm{P_l f}{L^q}\leq C_n 2^{l n (1/p-1/q)}\norm{P_l f}{L^p}.
$$
for $1\leq p<q\leq \infty$. \\
The (uncentered) Hardy-Littlewood maximal function is
$$
Mf(x)=\sup_{Q\supset x} \f{1}{|Q|} \int_Q |f(y)| dy.
$$
It is well-known that $M:L^p\to L^p$ for all $1<p\leq \infty$ and is of 
weak type $(1,1)$. It is also convenient to use the pointwise bound 
\begin{equation}
\label{eq:30}
\sup_{t>0} t^{-n} |f*\Phi(t^{-1} \cdot)(x)|\leq C \|\Phi\|_{L^1} Mf(x),
\end{equation}
for a radially dominated function $\Phi$. 
For integer values of $s$, we may define $W^{p,s}$ to be 
the Sobolev space with $s$ derivatives in $L^p$, $1\leq p\leq \infty$, with the corresponding norm  
$$
\norm{f}{W^{p,s}}:=\suml_{|\al|\leq s} \norm{D^\al_x f}{L^p}. 
$$
Equivalently, and for noninteger values of $s$, define 
$$
\norm{f}{W^{p,s}}:=\norm{f}{L^p}+ 
\|(\suml_{l=0}^\infty  2^{2l s}|P_l f|^2)^{1/2}\|_{L^p}
$$
and its  homogeneous analogue 
$$
\norm{f}{\dot{W}^{p,s}}:=
\|(\suml_{l=-\infty}^\infty  2^{2l s}|P_l f|^2)^{1/2}\|_{L^p}.
$$
Note $W^{p,s}=\dot{W}^{p,s}\cap L^p$.  \\
The  (homogeneous) Besov spaces $\dot{B}_{p,q}^s$, 
which scale like $\dot{W}^{p,s}$, are  defined as follows 
$$
\norm{f}{\dot{B}_{p,q}^s}:= 
(\suml_{l\in \cz}  2^{l sq }\|P_l f\|_{L^p}^q)^{1/q}.
$$
The Triebel-Lizorkin spaces are defined via 
$$
\norm{f}{\dot{F}_{p,q}^s}:= 
\|(\suml_{l\in \cz} 2^{l sq }|P_l f|^q)^{1/q}\|_{L^p}.
$$

\subsection{Fourier analysis on $\sn$} In this section, we define the Sobolev and Besov spaces for functions $q$ defined on $\sn$.  For that, the standard approach is to fix the basis of the spherical harmonics and define the Littlewood-Paley operators by projecting over the corresponding set of the harmonics within the fixed frequency.  

Introduce the angular differentiation operators $\Om_{i j}=x_j \p_i - x_i \p_j$.  It is well-known  that $\{\Om_{ i j}\}_{i\neq j}$ generate the algebra of all differential operators, acting on $C^\infty(\sn)$. The spherical Laplacian is defined via 
$$
\De_{sph}=\suml_{i< j} \Om_{i j}^2.
$$
The spherical harmonics $\{Y^n_{l,k}\}_{k\in A^n_l}$ 
 are eigenfunctions of $\De_{sph}$,  so that $\De_{sph} Y^n_{l,k}=-l(n-2+l) Y^n_{l,k}$, where $l\geq 0$, $k$ varies in a finite set $A_l^n$.  
  An equivalent way to define them is to take all the 
homogeneous of degree $l$ polynomials that are solutions to 
 \begin{equation}
 \label{eq:45}
 (\p_r^2+ r^{-1} \p_r + r^{-2} \De_{sph})Y^l=0
 \end{equation}
 Iy turns out that \eqref{eq:45} 
has $\left(\begin{array}{c} n+l-1 \\ l\end{array}\right)$ 
linearly independent solutions 
 $\{Y^n_{l,k}\}_{k\in A^n_l}$.  
Another  important property of the family  $\{Y^n_{l,k}\}$ 
 is that it forms 
  an orthonormal  basis for $L^2(\sn)$.  
  
 Let $f:\sn\to \cc$ be a smooth function.  One can then define 
the expansion in spherical harmonics in the usual way 
$$
f(\theta)=\suml_{l,k\in A_l^n} c_{l,k}^n Y^n_{l,k}(\theta), 
$$
where $c_{l,k}^n=\dpr{f}{Y^n_{l,k}}_{L^2(\sn)}$. 
 The Littlewood-Paley operators may be  defined via 
 $$
 P_m^{\xi/|\xi|} f = \suml_{l,k\in A_l^n} c_{l,k}^n\vp(2^{-m} l)  
Y^n_{l,k}(\theta), 
 $$
 and there is the equivalence for all\footnote{The constant 
of equivalence here depends only on $p$ and the cutoff 
function $\vp$.}  $1<p<\infty$, due to Strichartz \cite{Str1972}
 $$
 \norm{f}{L^p(\sn)}\sim 
\|(\sum_{k=-\infty}^\infty |P^{\xi/|\xi|}_k f|^2 )^{1/2}\|_{L^p(\sn)}.
 $$
 As a simple consequence, one has for all $1<p<\infty$, 
 $\|P_k^{\xi/|\xi|}\|_{L^p\to L^p}\leq C_{n,p}$. 
Such an inequality actually extends 
(as in the case of $\rn$) to the endpoint cases $p=1, 
p=\infty$, see \cite{SteinWeiss}. 
 One can also define the Sobolev spaces $W^{p,s}(\sn):1<p<\infty$ via 
 \begin{eqnarray*}
 & & \norm{f}{W^{p,s}(\sn)}=\suml_{|\al|\leq s}\norm{ \Om^\al f}{L^p(\sn)} \\
 & & \norm{f}{\dot{W}^{p,s}(\sn)}= 
\|(\suml_{k=-\infty}^\infty2^{2ks} |P^{\xi/|\xi|}_k f|^2 )^{1/2}\|_{L^p(\sn)}.
 \end{eqnarray*}
 where $\Om=\sqrt{-\De_{sph}}$. These last two formulas 
  give equivalent definitions for the case of integer $s$.  The (homogeneous) 
  Besov spaces are defined in the usual manner as follows 
 \begin{eqnarray*}
& &  \norm{f}{\dot{B}^s_{p,q}(\sn)}=
(\suml_{k=-\infty}^\infty 2^{qks} \|P^{\xi/|\xi|}_k f\|_{L^p(\sn)}^q)^{1/q}.
\end{eqnarray*}
It is worth mentioning at this point that a variant of Bernstein 
inequality holds\footnote{The proof is simply 
 that there are  $\sim N^{n-1}$ spherical harmonics at frequency 
$N$, just as for the Bernstein inequality one uses that the volume of 
$\{\xi\in\rn:|\xi|\sim N\}$ is $\sim N^n$.} in this context, 
see \cite{Sterbenz}, p. 201. This together with the Littlewood-Paley theory outlined above implies Sobolev embedding for $L^p(\sn)$ spaces. 
For future reference,  let us record  
this estimate 
\begin{equation}
\label{eq:bern}
\|f\|_{L^q(\sn)}\leq C_{p,q,n} 
\norm{\Omega^{(n-1)(1/p-1/q)}f}{L^p(\sn)}, 
\end{equation}
which holds whenever for $1<p< q< \infty$.

If $q$ is a finite sum of harmonics, it is actually an analytic function 
  (it is a fact a restriction of a polynomial to the unit sphere) and one may write 
\begin{equation}
\label{eq:36}
q(\theta)=\suml_{\al>0} \f{D_\xi^\al q(\theta_0)}{\al!}(\theta-\theta_0)^\al.
\end{equation}
Here, $D_\xi^\al q(\theta_0)$ should be understood as taking $\al$ derivatives of 
the corresponding homogeneous polynomial and evaluating at $\theta_0$. 
The following lemma is standard, but since we need a specific dependence of our 
estimates  upon the parameter $\al$, we state it here for completeness. 
\begin{lemma}
\label{le:os}
Let $q:\sn\to \cc$ and $q_m=P^{\xi/|\xi|}_m q$. Then, there is a constant $C_n$, so that 
for every $1\leq p\leq \infty$, there is a the estimate 
$$
\|D^\al_\xi q_m\|_{L^p(\sn)}\leq C_n^{|\al|} 2^{m |\al|} \norm{q_m}{L^p(\sn)}.
$$
\end{lemma}
The proof of Lemma \ref{le:os} is standard. One way to proceed is to note that if we
extend the function $q_m$ off $\sn$ to some annulus, say  via 
$Q_m(\xi)=\vp(|\xi|) q_m(\xi/|\xi|)$, and 
then 
$$
 \|D^\al_\xi q_m\|_{L^p(\sn)}\lesssim \|D^\al_\xi Q_m \|_{L^p(\rn)}
$$

\section{$L^p$ estimates for PDO with rough symbols}
We start with the $L^2$ estimate to illustrate the main ideas in the proof. 
\subsection{$L^2$ estimates: Proof of Theorem \ref{theo:1}}
Our first remark is that we will for convenience consider only 
real-valued symbols $\si$, since of course the general 
case follows from splitting into a real and imaginary part. 

To show $L^2$ estimates for $T_\si$, it is equivalent to show $L^2$ 
estimates for the adjoint operator, which takes the 
form\footnote{ There is the small technical problem 
that the $\xi$ integral does not converge 
absolutely. 
This can be resolved by judicious placement of cutoffs 
$\chi(\xi/N)$, after which, one may subsume that 
part in $\hat{f}(\xi)$. In  the end 
we let $N\to \infty$ and all the estimates will be 
independent of the cutoff constant $N$.} 
$$
T_\si^* g(x)=\int e^{2\pi i \xi\cdot x} (\int e^{-2\pi i \xi\cdot y} 
[g(y) \si(y, \xi)] dy )d\xi.
$$
Our next task is to decompose $T_\si^* g$ and we start by 
taking a Littlewood-Paley partition of unity in the $\xi$ variable for $g$. We have 
$$
T_\si^* g (x)=\suml_{l\in \cz} \int e^{2\pi i \xi\cdot x} 
(\int e^{-2\pi i \xi\cdot y} 
[ g(y) P_l^\xi \si(y, \xi)] dy )d\xi
$$ 
Now that the function $g$ is frequency localized at frequency $2^l$, 
we  introduce further decomposition in the  
$\xi$ integration. 

For  the $L^2$ estimates, because of the orthogonality, 
we only need  rough partitions, so for each fixed $l$, 
take a tiling of $\rn$ composed of cubes $\{Q\}$ 
with diameter $2^{-l}$. Denote the  characteristic functions of $Q$ by $\chi_Q$. We have 
$$
T_\si^* g (x)=\suml_{l\in \cz} \suml_{Q:d(Q)=2^{-l}} 
\int e^{2\pi i \xi\cdot x} \chi_Q(\xi)
(\int e^{-2\pi i \xi\cdot y} 
[ g(y) P_l^\xi \si(y, \xi)] dy )d\xi
$$
The main point of our next decompositions is that the 
function $P_l^\xi \si$ is essentially constant in $\xi$ 
over any fixed cube $Q$. We exploit that by observing  that 
$\xi\to P_l^\xi\si(x, \xi)$ is an entire  function and there is 
 the expansion 
$$
P_l^\xi \si(y, \xi) \chi_Q(\xi)=[P_l^\xi \si(y, \xi_Q)+\suml_{\al: |\al|>0}^\infty
\f{D_\xi^\al P_l^\xi \si(y, \xi_Q)}{\al!} (\xi-\xi_Q)^\al]\chi_Q(\xi)
$$
for any fixed $y$ and for any $\xi_Q\in Q$. Note that 
$D_\xi^\al P_l^\xi \si(y, \xi_Q)\sim 2^{l|\al|} P_l^\xi \si(y, \xi_Q)$ and 
$ |(\xi-\xi_Q)^\al|\lesssim  2^{-l|\al|}$, by support consideration 
(recall $d(Q)=2^{-l}$). On a heuristic level, by the presence of  $\al!$, 
one should think that  the series 
above behave like $P_l^\xi \si(y, \xi_Q)$ plus exponential tail. 

Going back to $D_\xi^\al P_l^\xi$, as  we have mentioned in Section \ref{sec:prelim}, we can write 
$D_\xi^\al P_l^\xi=2^{l|\al|} P_{l, \al}^\xi$, where 
$P_{l, \al}^\xi$ is given by the multiplier $\vp(2^{-l}\xi) (2^{-l} \xi)^\al$. 
It is clear that 
$\|P_{l, \al}^\xi f\|_{L^2(\rn)}\leq C_n^{|\al|}\|P_l^\xi f\|_{L^2(\rn)}$. 

Thus, we have arrived at 
\begin{eqnarray*}
& & T_\si^* g (x)=\suml_{|\al|\geq 0} \suml_{l\in \cz} \suml_{Q:d(Q)=2^{-l}} 
\int e^{2\pi i \xi\cdot x} \chi_Q(\xi) 
\f{2^{l|\al|}(\xi-\xi_Q)^\al}{\al!} \times \\
& & \times 
(\int e^{-2\pi i \xi\cdot y} 
[ g(y) P_{l,\al}^\xi \si(y, \xi_Q)] dy )d\xi= 
\suml_{l,\al}(\al!)^{-1} \suml_{l, Q:d(Q)=2^{-l}}  P_{Q, l, \al} [ g(\cdot) 
P_{l,\al}^\xi \si(\cdot, \xi_Q)], 
\end{eqnarray*}
where $P_{Q, l, \al}$ acts via $\widehat{P_{Q, l, \al} f}(\xi)= \chi_Q(\xi) 
2^{l|\al|}(\xi-\xi_Q)^\al \hat{f}(\xi)$. Note 
$$
\|P_{Q, l, \al}\|_{L^2\to L^2}=\sup_\xi | \chi_Q(\xi) 
2^{l|\al|}(\xi-\xi_Q)^\al|\leq 1. 
$$

For fixed $l, \al$, take $L^2$ norm. Using the orthogonality of 
$P_{Q, l, \al}$ and its boundedness on $L^2$, we obtain 
\begin{eqnarray*}
& & \|\suml_{Q:d(Q)=2^{-l}}  P_{Q, l, \al} [ g(\cdot) 
P_{l,\al}^\xi \si(\cdot, \xi_Q)]\|_{L^2}^2= 
\suml_{Q:d(Q)=2^{-l}} \|P_{Q, l, \al} [ g(\cdot) 
P_{l,\al}^\xi \si(\cdot, \xi_Q)]\|_{L^2}^2\leq  \\
& & \leq \suml_{Q:d(Q)=2^{-l}} \|g(\cdot) 
P_{l,\al}^\xi \si(\cdot, \xi_Q)]\|_{L^2}^2=\int |g(y)|^2 
(\suml_{Q} |P_{l,\al}^\xi \si(y, \xi_Q)|^2) dy . 
\end{eqnarray*}
We now again use   $P_{l,\al}^\xi \si(y, \xi_Q)\sim P_{l,\al}^\xi \si(y, \eta)$ for any
$\eta\in Q$, this time to estimate the contribution of $\sum_{Q} |P_{l,\al}^\xi \si(y, \xi_Q)|^2$. 
This is done as follows. Expand 
\begin{equation}
\label{eq:859}
P_{l, \al} ^\xi \si(y, \xi_Q)= \suml_{\beta: |\beta|\geq 0}^\infty
\f{D_\eta^\be P_{l, \al}^\eta \si(y, \eta)}{\be!} (\xi_Q-\eta)^\be,
\end{equation}
to be used for  $\eta\in Q$. Thus, if we average over $Q$, 
\begin{eqnarray*}
& & |P_{l, \al} ^\xi \si(y, \xi_Q)|=\left(|Q|^{-1} \int_Q |\suml_{\beta: |\beta|\geq 0}^\infty
\f{D_\xi^\be P_{l, \al}^\xi \si(y, \eta)}{\be!} (\xi_Q-\eta)^\be|^2 d\eta
\right)^{1/2}\leq  \\
& &\leq |Q|^{-1/2} \suml_{\beta: |\beta|\geq 0}^\infty \f{C_n^{|\be|} 
2^{-l|\be|}}{\be!}  \left(\intl_Q | D_\xi^\be P_{l, \al}^\xi \si(y, \eta)|^2 d\eta\right)^{1/2}.
\end{eqnarray*}
and so (recalling $|Q|\sim 2^{-l n}$)
\begin{eqnarray*}
& &(\suml_{Q} |P_{l,\al}^\xi \si(y, \xi_Q)|^2)^{1/2}\leq C 2^{l n/2} \suml_{\be} \f{C_n^{|\be|} 
2^{-l|\be|}}{\be!} \|D_\xi^\be P_{l, \al}^\xi \si(y, \cdot)\|_{L^2} \leq \\
& & \leq 
C 2^{l n/2} \|P_{l, \al}^\xi \si(y, \cdot)\|_{L^2}.
\end{eqnarray*}
Thus, 
\begin{eqnarray*}
& & \|T_\si^* g\|_{L^2} \lesssim  \suml_{l, \al} 2^{l n/2} (\al!)^{-1} \left(\int |g(y)|^2   
 \|P_{l, \al}^\xi \si(y, \cdot)\|_{L^2}^2 dy\right)^{1/2}\lesssim \\
 & &\lesssim 
 \suml_{l, \al} 2^{l n/2} (\al!)^{-1} \|g\|_{L^2} \supl_y \|P_{l, \al}^\xi \si(y, \cdot)\|_{L^2}. 
\end{eqnarray*}
Furthermore, 
$$
\sup_y \|P_{l,\al}^\xi \si(y,\cdot)\|_{L^2}\leq C_n^{|\al|} 
\sup_y \|P_{l}^\xi \si(y,\cdot)\|_{L^2}.
$$
Put everything together 
\begin{eqnarray*}
& & \|T_\si^* g\|_{L^2}\leq  C_n 
\norm{g}{L^2}\suml_{\al} (\al!)^{-1} C_n^{|\al|} 
\sum_l 2^{ln/2} 
 \sup_y \|P_{l}^\xi \si(y,\cdot)\|_{L^2} 
\leq \\
& & \leq D_n \norm{g}{L^2} \sum_l 2^{ln/2} 
 \sup_y \|P_{l}^\xi \si(y,\cdot)\|_{L^2}, 
\end{eqnarray*}
as desired. 
\subsection{$L^p$ estimates: $2<p\leq \infty$}
The result anounced in Theorem \ref{theo:3} 
follows by interpolation between the $L^2$ estimate just 
proved and the boundedness of $T_\si:L^\infty\to L^\infty$, 
which we need to show next. 

We do  that by showing that the adjoint $T_\si^*:L^1\to L^1$.  
This is relatively  easy, since one can reduce to showing that 
\begin{equation}
\label{eq:31}
\sup_{a\in\rn} \norm{T_\si^* \de(a)}{L^1(\rn)}<C\sum_l 2^{ln/2} 
 \sup_a \|P_{l}^\xi \si(a,\cdot)\|_{L^2}.
\end{equation}

This is standard,  since one can embed $L^1$  into the space 
of all Borel measures $M(\rn)$. The next observation is that 
by Krein-Milman's theorem, the convex combinations of 
the set of Dirac masses $\{\de_a: a\in \rn\}$  are weak* dense 
in the unit ball of $M(\rn)$. 

Fix $a\in \rn$. We have 
\begin{eqnarray*}
& & T_\si^* \de_a(x)= \int e^{2\pi i \xi\cdot(x-a)} \si(a, \xi) d\xi=\cf_\xi(\si(a, \cdot))(a-x).
\end{eqnarray*}
where $\cf_\xi$ signifies the Fourier transform in the $\xi$ 
variable. Denote $g_a(z)=\cf_\xi(\si(a, \cdot))(z)$. 
We have by Cauchy-Schwartz and the Plancherel's theorem 
\begin{eqnarray*}
& & \|T_\si^* \de_a\|_{L^1}=\int |g_a(a-x)| dx=
\suml_{l\in \cz} \int_{|x-a|\sim 2^l} |g_a(a-x)| dx\leq  \\
& &\leq C_n 
\suml_{l\in \cz} 2^{ln/2} 
(\int_{z: |z|\sim 2^l} |g(z)|^2 dz)^{1/2}\leq  C_n 
\suml_{l\in \cz} 2^{ln/2} \|P_l^\xi \si(a, \cdot)\|_{L^2}
\end{eqnarray*}
This is \eqref{eq:31}, whence \eqref{eq:22}.  

\subsection{$L^p$ estimates: $1<p\leq 2$}
We take slightly different approach than in the case of $L^2$ estimates. 
Namely, we will show  that $T_\si$ 
is of weak type $(1,1)$ operator, whence, by interpolation with the 
$L^2$ estimate, one gets the full range $1<p\leq 2$. 
Note that the $L^2$ estimate comes with 
$$
\norm{T_\si}{L^2\to L^2}\leq C_n \suml_l 
2^{ln/2} \sup_y \|P_l^\xi \si(y, \cdot)\|_{L^2}
\leq C_n \suml_l  2^{l n} \sup_y \|P_l^\xi \si(y, \cdot)\|_{L^1}, 
$$
where in the last inequality, we have used the Bernstein inequality \\
$\|P_l^\xi \si(a, \cdot)\|_{L^2}\leq C_n 2^{l n/2} 
\|P_l^\xi \si(a, \cdot)\|_{L^1}$. Thus, it remains to show weak type $(1,1)$ bounds for $T_\si$. 
We proceed by performing a decomposition for $T_\si$, 
 inspired by the $L^2$ bounds. Our goal is to show the pointwise estimate 
\begin{equation}
\label{eq:32}
|T_\si f(x)|\leq C_n (\suml_{l\in \cz} 
2^{l n}\sup_y  \|P_l^\xi \si(y, \cdot)\|_{L^1}) Mf(x),
\end{equation}
which implies the desired weak type bounds since $M:L^1\to L^{1, \infty}$. 

To achieve that, we have to be a bit more careful than in the $L^2$ case, 
since the rough cutoffs in the $\xi$ variable will be insuficient to 
show \eqref{eq:32}. 

For any integer $l$, introduce smooth partition of unity, which is 
adapted to 
the cover $\rn=\cup_{Q: diam(Q)= 2^{-l}}Q$, that is a family of 
functions $\{h_{l,Q}\}$, with 
$supp\  h_{l, Q}\subset Q^*$ and 
$ |D^\al h_{l, Q}(\xi)|\leq C_\al 2^{l|\al|}$ for 
every multiindex $\al$.  Choose and fix a family of arbitrary points $\xi_Q\in Q$. 
By rescaling, 
one can choose $h_{l, Q}:=\psi_{l,Q}(2^{l}(\xi-\xi_Q))$, 
where $supp\ \psi_{l, Q}\subset \{|\xi|<2\}$ and  $ |D^\be_{\eta}
\psi_{l, Q}(\eta)|\leq C_\be$ and 
$$
\suml_{Q}\psi_{l, Q}(2^l(\xi-\xi_Q))=1
$$ 
Write 
\begin{eqnarray*}
& & T_\si f(x)=\suml_{l\in \cz} 
\intl_{\rn} P_l^\xi \si(x,\xi) e^{2\pi i \xi x} \hat{f}(\xi) d\xi= \\
& &= 
\suml_{l\in \cz}\suml_Q \intl_{\rn} P_l^\xi \si(x,\xi)
\psi_{l,Q}(2^l(\xi-\xi_Q)) 
 e^{2\pi i \xi x} \hat{f}(\xi) d\xi
\end{eqnarray*}
We now expand the $P_l^\xi \si(x,\xi)$ around $\xi_Q$. We have 
$$
P_l^\xi \si(x, \xi) \psi_{l,Q}(2^l(\xi-\xi_Q))= 
(\sum_{\al: |\al|\geq 0}^\infty
\f{D_\xi^\al P_l^\xi \si(x, \xi_Q)}{\al!} (\xi-\xi_Q)^\al) 
\psi_{l,Q}(2^l(\xi-\xi_Q))
$$
Plugging that in the formula for $T_\si f$ yields 
\begin{eqnarray*}
& & T_\si f(x)=\suml_{l, \al} (\al!)^{-1}   \suml_Q 2^{-l|\al|}
D_\xi^\al P_l^\xi \si(x, \xi_Q)  
Z_{l, Q}^{\al} f(x). 
\end{eqnarray*}
where $\widehat{Z_{l, Q}^{\al} f}(\xi)=
\psi_{l,Q} (2^l(\xi-\xi_Q)) (2^l(\xi-\xi_Q))^{\al}
 \hat{f}(\xi)=\psi^\al_{l,Q} (2^l(\xi-\xi_Q))\hat{f}(\xi)$, i.e. \\
$\psi^\al_{l,Q}(z)=\psi_{l,Q}(z) z^\al$. 
By \eqref{eq:30}, we get 
$$
|Z_{l, Q}^{\al} f(x)|\leq C_n \|\widehat{\psi^\al_{l,Q}}\|_{L^1} M f(x).
$$
By the elementary properties of the Fourier transform 
\begin{eqnarray*}
& &\|\widehat{\psi^\al_{l,Q}}\|_{L^1} = C_n^{|\al|}  \int |D^\al_\eta[
 \widehat{\psi_{l,Q}}(\eta)]| d\eta\leq C_n^{|\al|}  \sum_{k=-\infty}^\infty  2^{k |\al|} 
\int |P_k^\eta[ \widehat{\psi_{l,Q}}(\eta)]| d\eta
\end{eqnarray*}
But by support considerations, 
 $P_k^\eta [\widehat{\psi_{l,Q}}]=0$ if $k>3$. Also  since  $P_k:L^1\to L^1$, we get 
$$
\|\widehat{\psi^\al_{l,Q}}\|_{L^1}\leq C_n^{|\al|} 
\|P_k[\widehat{\psi_{l,Q}}]\|_{L^1}\leq  C_n^{|\al|} 
\|\widehat{\psi_{l,Q}}\|_{L^1}\leq 
C_n^{|\al|}.
$$
Thus, it remains to show for every $x$  and for {\it any}
$\{\xi_Q\}, \xi_Q\in Q$
\begin{equation}
\label{eq:35} 
\suml_\al \f{2^{-l|\al|}}{\al!} \suml_Q 
|D_\xi^\al P_l^\xi \si(x, \xi_Q)|  \leq C_n 
2^{ln} \sup_y \|P_l^\xi \si(y, \cdot)\|_{L^1(\rn)}. 
\end{equation}
This is done similar to the $L^2$ case. By  \eqref{eq:859} and by averaging over 
the corresponding $Q$
\begin{eqnarray*}
& & \suml_\al \f{2^{-l|\al|}}{\al!} \suml_Q 
|D_\eta^\al P_l^\eta \si(x, \xi_Q)|\leq \suml_{\al, \be} \f{2^{-l|\al|}}{\al!\be!}  
 \suml_Q |Q|^{-1} \int_Q  |D_\eta^{\al+\be} P_l^\eta \si(x, \eta)(\eta-\xi_Q)^\be|d\eta\\
 & & \lesssim 2^{ln}  \suml_{\al, \be} \f{2^{-l|\al|}}{\al!\be!} 
 \int |D_\eta^{\al+\be} P_l^\eta \si(x, \eta)(\eta-\xi_Q)^\be|d\eta\lesssim 2^{ln} \|P_l^\eta \si(x,
 \cdot)\|_{L^1} 
\end{eqnarray*}

\section{$L^p$ estimates for homogeneous of degree zero symbols}
We start with the $L^2$ estimate, since it is very similar to the corresponding 
estimate \eqref{eq:7} and contains the main ideas for the $L^p$ estimate. 

\subsection{$L^2$ estimates for homogeneous of degree zero symbols}
Consider $T_\si^*$ and introduce the Littlewood-Paley partition of unity $P_l^{\xi/|\xi|}$. We have 
$$
T_\si^* g (x)=\suml_{l=0}^\infty \int e^{2\pi i \xi\cdot x} 
(\int e^{-2\pi i \xi\cdot y} 
[ g(y) P_l^{\xi/|\xi|}  q(y, \xi/|\xi|)] dy )d\xi
$$ 
For every $l\geq 0$, introduce a partition of unity on $\sn$, say $\{K\}$, which consists of disjoint 
sets of diameter comparable to $2^{-l}$. One may form $\{ K\}$ by introducing a $2^{-l}$ net on $\sn$, say $\xi^m_{l}$,  form  the conic sets  $H_m^l=\{\xi\in\rn:  \: |\xi/|\xi|-\xi^m_l|\leq 2^{-l}\}$ and construct \\ $K^l_m=H_m^l\setminus \cup_{j=0}^{m-1} H_{j}^l$.   We have 
\begin{equation}
\label{eq:400}
T_\si^* g (x)=\suml_{l=0}^\infty  \suml_m \int e^{2\pi i \xi\cdot x} \chi_{K^l_m}(\xi)
(\int e^{-2\pi i \xi\cdot y} 
[ g(y) P_l^{\xi/|\xi|} q (y, \xi/|\xi|)] dy )d\xi
\end{equation}
Now, that the symbol is frequency localized around frequencies $\sim 2^l$ 
 and the sets $K^l_m\cap \sn$ have diameters less 
than $2^{-l}$,  we expand $q(y, \xi/|\xi|)$ around an {\it arbitrary point} $\theta_l^m\in K_l^m$.
According to \eqref{eq:36}, we have for all $\xi\in K^l_m$, 
$$
q(y, \xi/|\xi|) 
=\suml_{\al\geq 0} \f{D_\xi^\al q(y,\theta_l^m)}{\al!}(\xi/|\xi|-\theta_l^m)^\al.
$$
Entering this new expression in \eqref{eq:400} yields 
\begin{eqnarray*}
& & 
T_\si^* g (x)=\suml_{l=0}^\infty  \suml_m \suml_{\al} (\al!)^{-1} 
\int e^{2\pi i \xi\cdot x} \chi_{K^l_m}(\xi) (\xi/|\xi|-\theta_l^m)^\al \times \\
& & \times 
(\int e^{-2\pi i \xi\cdot y} 
[ g(y)  D_\xi^\al P_l^{\xi/|\xi|} q (y, \theta_l^m)] dy )d\xi= \\
& & =
\suml_{l=0}^\infty  \suml_m \suml_{\al} (\al!)^{-1}  Z_{l,m}^\al [g(\cdot) 2^{-l|\al|} 
 D_\xi^\al P_l^{\xi/|\xi|} q (\cdot , \theta_l^m)] ,
\end{eqnarray*}
where $Z_{l,m}^\al$ is given by the multiplier $\chi_{K^l_m}(\xi) 2^{l |\al|} 
(\xi/|\xi|-\theta_l^m)^\al$.  Note the disjoint support of the multipliers 
$\{Z_{l,m}^\al\}_m$ and $\|Z_{l, m}^\al\|_{L^2\to L^2}=\sup_\xi|\chi_{K^l_m}(\xi) 2^{l |\al|} 
(\xi/|\xi|-\theta_l^m)^\al| \leq 4^{|\al|}$. \\
Take $L^2$ norm of $T_\si^* g$.   
\begin{eqnarray*}
& & \norm{T_\si^* g}{L^2(\rn)}\lesssim 
\suml_{l=0}^\infty  \suml_{\al} (\al!)^{-1}\left(\suml_m 
\norm{Z_{l, m}^\al   [g(\cdot) 2^{-l|\al|} 
 D_\xi^\al P_l^{\xi/|\xi|} q (\cdot , \theta_l^m) ]}{L^2}^2 \right)^{1/2}\leq \\
& &\leq  4^{|\al|} \suml_{l=0}^\infty \suml_{\al} (\al!)^{-1} \left(\suml_m \norm{g(\cdot) 2^{-l|\al|} 
 D_\xi^\al P_l^{\xi/|\xi|} q (\cdot , \theta_l^m) }{L^2}^2\right)^{1/2}.
\end{eqnarray*}
We proceed to further bound the expression in the  $m$ sum. Since
$$
\suml_m \norm{g(\cdot) 2^{- l|\al|} 
 D_\xi^\al P_l^{\xi/|\xi|} q (\cdot , \theta_l^m)}{L^2}^2=2^{-2 l|\al|} \int |g(y)|^2 \left(\suml_m 
| D_\xi^\al P_l^{\xi/|\xi|} q (y , \theta_l^m)|^2\right) dy, 
$$
matters reduce to a good estimate for $\suml_m 
| D_\xi^\al P_l^{\xi/|\xi|} q (y , \theta_l^m)|^2$. 
We proceed as before. By \eqref{eq:36}, we get for all $\eta\in K^l_m\cap \sn$, 
\begin{equation}
\label{eq:n3}
  D_\xi^\al P_l^{\xi/|\xi|} q(y, \theta_l^m) 
=\suml_{\be\geq 0} \f{D_\xi^{\al+\be} P_l^{\eta/|\eta|}q(y,\eta)}{\be!}(\theta_l^m-\eta)^\be.
\end{equation}
Averaging over $K^l_m\cap \sn$ and taking into account $|K^l_m\cap \sn|\sim 2^{l(n-1)}$  yields 
\begin{eqnarray*}
& & (\suml_m 
| D_\xi^\al P_l^{\xi/|\xi|} q (y , \theta_l^m)|^2)^{1/2} \lesssim \\
& & \lesssim \suml_{\be} \f{2^{-l |\be|}}{\be!}
(\suml_m 
|K^l_m\cap \sn|^{-1} \intl_{K^l_m\cap \sn} |D_\xi^{\al+\be} P_l^{\eta/|\eta|}q(y,\eta)|^2
d\eta)^{1/2}\\
& &\lesssim 2^{l(n-1)/2} \suml_{\be} \f{2^{-l |\be|}}{\be!} 
\|D_\xi^{\al+\be} P_l^{\eta/|\eta|}q(y,\cdot)\|_{L^2} \\
& &\lesssim 2^{l[(n-1)/2+|\al|]} 
\suml_{\be} \f{C_n^{|\al|+|\be|}}{\be!} 
\|P_l^{\eta/|\eta|}q(y,\cdot)\|_{L^2(\sn)} \\ 
& &\leq C_n^{|\al|}2^{l[|\al|+(n-1)/2]}
\|P_l^{\eta/|\eta|}q(y,\cdot)\|_{L^2(\sn)}.
\end{eqnarray*}
Putting this back into the estimate for $\|T_\si^* g\|_{L^2(\rn)}$ implies 
$$
\norm{T_\si^* g}{L^2(\rn)}\lesssim \|g\|_{L^2} \suml_l 2^{l(n-1)/2}  \supl_y 
\|P_l^{\eta/|\eta|}q(y,\cdot)\|_{L^2(\sn)}
$$
as desired. 
\subsection{$L^p$ estimates for homogeneous of degree zero multipliers}
Fix $p: 2\leq p<\leq \infty$. 
To verify the estimate $\|T\|_{B^0_{p,1}\to L^p}$, it will suffice to fix $k$ and show 
\begin{equation}
\label{eq:01}
\|T(P_k f)\|_{L^p}\leq C \|f\|_{L^p}.
\end{equation}
Furthermore, by the scale invariance of the quantity  $\sum_l 2^{l(n-1)} 
\sup_y\|P_l^{\xi/|\xi|} q(y, \cdot)\|_{L^1(\sn)}$  this is equivalent to  verifying \eqref{eq:01} 
only for  $k=0$. 
That is, it suffices to establish 
the $L^p, p\geq 2$ boundedness of the  operator   
$$
G f(x)=\intl_{\rn} q(x,\xi/|\xi|) e^{2\pi i \xi x} \vp(|\xi|) \hat{f}(\xi) d\xi.
$$
provided the multiplier  $m$ satisfies 
$\suml_l 2^{l(n-1)} \sup_y\|P_l^{\xi/|\xi|} q(y, \cdot)\|_{L^1(\sn)}<\infty$. 

Next, we make the angular decomposition as in the case of the $L^2$ estimates for the adjoint 
operator $G^*$.However, this time we will have to be more careful and 
instead of the rough cutoffs $\chi_{K^l_m}$, we shall use a 
smoothed out versions of them. Fix $l$.  Choose and fix a 
 $2^{-l}$ net  $\theta_m^l\in K_m^l\cap \sn$, so that the family 
$\{\theta\in \sn: 
|\theta_m^l-\theta|\leq 2^{-l}\}_m$ has the finite intersection
property. 
Introduce a family of functions $\vp_{l,m}:\rn\to [0,1]$, so that for every $\xi\in \rn$, 
\begin{eqnarray}
\label{eq:fun}
& & \suml_m \vp_{l,m} (2^l(\xi/|\xi|-\theta_m^l))=1 \\ 
\nonumber
& & \sup_\eta |D^\be_\eta \vp_{l,m} (\eta)|\leq C_\be.
\end{eqnarray}
In other words, the family of functions $\{\vp_{l,m}\}$ provides a 
smooth partition of unity, subordinated to the cover $\{K_m^l\}$. 

As before, write 
$$
G^* g(x)=\suml_{l\geq 0} \intl_{\rn} e^{2\pi i \xi x}  \vp(|\xi|)  \int  e^{-2\pi i \xi y} [g(y) 
P_l^{\xi/|\xi|} q(y, \xi/|\xi|)]dy  d\xi.
$$
Inserting the partition of unity discussed above   into the ($l^{th}$ term of the) 
last formula for $G^*$ yields 
$$
G^* g (x)=\suml_{l\geq 0} \suml_m  \int  e^{2\pi i \xi (x-y)} 
\vp_{l,m} (2^l(\xi/|\xi|-\theta_m^l)) \vp(|\xi|)  [g(y) 
P_l^{\xi/|\xi|} q(y, \xi/|\xi|)]dy d\xi.  
$$
Following the same strategy as before, we expand $q(y, \xi/|\xi|)$ around  $\theta_m^l\in K_m^l$.
According to \eqref{eq:36}, we have 
$$
P_l^{\xi/|\xi|} q(y, \xi/|\xi|) 
=\suml_{\al\geq 0} \f{D_\xi^\al P_l^{\xi/|\xi|} 
q(y, \theta_m^l)  }{\al!}(\xi/|\xi|-\theta_m^l)^\al.
$$
Of course,  the last formula is useful only when 
$|\xi/|\xi|-\theta_m^l|\lesssim  2^{-l}$, in particular on the 
support of $\vp_{l,m} (2^l(\xi/|\xi|-\xi_m^l)) $. 
This gives us the representation 
\begin{eqnarray*}
& & G^* g =\suml_{l\geq 0}\suml_m   \suml_{|\al|\geq 0} (\al!)^{-1}
\int   e^{2\pi i \xi x} 
\vp_{l,m} (2^l(\xi/|\xi|-\theta_m^l)) (\xi/|\xi|-\theta_m^l)^\al \vp(|\xi|) \times \\
& &\times \int e^{2\pi i \xi y} g(y) 
P_l^{\xi/|\xi|} D_\xi^\al P_l^{\xi/|\xi|}  q(y, \theta_m^l)dy d\xi= \\
& & =
\suml_{l\geq 0}\suml_m   \suml_{|\al|\geq 0}(\al!)^{-1} Z_{l,m}^\al 
[g(\cdot)  2^{-l|\al|} D_\xi^\al  P_l^{\xi/|\xi|}  q(\cdot, \theta_m^l)]
\end{eqnarray*}
where 
$$
\widehat{Z_{l,m}^\al f}(\xi)=  \vp_{l,m} (2^l(\xi/|\xi|-\theta_m^l))  2^{l|\al|} 
(\xi/|\xi|-\theta_m^l)^\al 
\vp(|\xi|) \hat{f}(\xi)=\vp_{l,m}^\al (\xi/|\xi|-\theta_m^l)\vp(|\xi|) \hat{f}(\xi).
$$
Taking $L^p$ norm of $G^* g$, we get 
\begin{eqnarray*}
& & \|G^* g\|_{L^p}\leq \suml_{l\geq 0} \suml_{|\al|\geq 0}(\al!)^{-1}\|\suml_m Z_{l,m}^\al 
[g(\cdot)  2^{-l|\al|} D_\xi^\al  P_l^{\xi/|\xi|}  q(\cdot, \theta_m^l)]\|_{L^p}
\end{eqnarray*}
Lemma \ref{le:sum} in the Appendix allows us to treat   expressions  of the type 
$\|\suml_m  Z_{l,m}^\al g_m^\al\|_{L^p}$. Indeed, according to \eqref{eq:n1}, we have 
\begin{eqnarray*}
& & \|\suml_m Z_{l,m}^\al 
[g(\cdot)  2^{-l|\al|} D_\xi^\al  P_l^{\xi/|\xi|} q(\cdot, \theta_m^l)]\|_{L^p}\lesssim  (\suml_m  
\|g(\cdot)  2^{-l|\al|} D_\xi^\al  P_l^{\xi/|\xi|}  q(\cdot, \theta_m^l)\|_{L^p}^p)^{1/p} \\
& &= 2^{-l|\al|} (\int |g(y)|^p (\suml_m  
   |D_\xi^\al  P_l^{\xi/|\xi|}  q(y, \theta_m^l)|^p) dy)^{1/p} 
\end{eqnarray*}
By virtue of  \eqref{eq:n3}, we get 
\begin{eqnarray*}
& & 
D_\xi^\al P_l^{\xi/|\xi|} q(y, \theta_l^m) 
=\suml_{\be\geq 0} \f{D_\xi^{\al+\be} P_l^{\eta/|\eta|}q(y,\eta)}{\be!}(\theta_l^m-\eta)^\be.
\end{eqnarray*}
whence by averaging\footnote{this step is identical to the one performed earlier for the 
$L^2$ bounds, except that now the $l^2$ sums are replaced by $l^p$ sums.} over $K^l_m\cap \sn$, 
\begin{eqnarray*}
 & & (\suml_m 
| D_\xi^\al P_l^{\xi/|\xi|} q (y , \theta_l^m)|^p)^{1/p} \lesssim \\
& & \lesssim \suml_{\be} \f{2^{-l |\be|}}{\be!}
(\suml_m 
|K^l_m\cap \sn|^{-1} \intl_{K^l_m\cap \sn} |D_\xi^{\al+\be} P_l^{\eta/|\eta|}q(y,\eta)|^p
d\eta)^{1/p}\\
& &\lesssim 2^{l(n-1)/p} \suml_{\be} \f{2^{-l |\be|}}{\be!} 
\|D_\xi^{\al+\be} P_l^{\eta/|\eta|}q(y,\cdot)\|_{L^p(\sn)} \\
& &\lesssim 2^{l[(n-1)/p+|\al|]} 
\suml_{\be} \f{C_n^{|\al|+|\be|}}{\be!} 
\|P_l^{\eta/|\eta|}q(y,\cdot)\|_{L^p(\sn)} \\ 
& &\leq C_n^{|\al|}2^{l[|\al|+(n-1)/p]}
\|P_l^{\eta/|\eta|}q(y,\cdot)\|_{L^p(\sn)}.
\end{eqnarray*}
All in all, 
\begin{eqnarray*}
& & \|G^* g\|_{L^p}\leq C_n \|g\|_{L^p}\suml_{l\geq 0} \suml_{|\al|\geq 0}(\al!)^{-1} 
 C_n^{|\al|}2^{l(n-1)/p}
\supl_y \|P_l^{\eta/|\eta|}q(y,\cdot)\|_{L^p(\sn)}\\ 
& &\leq C_n \|g\|_{L^p} \suml_l 
2^{l(n-1)/p}
\supl_y \|P_l^{\eta/|\eta|}q(y,\cdot)\|_{L^p(\sn)}.
\end{eqnarray*}
as desired. 

\section{Counterexamples}
\label{sec:counter}
\subsection{Theorem \ref{theo:1} is sharp.}
Given $p>2$, we will construct 
 an explicit symbol $\si(x, \xi)$, 
so that the corresponding PDO $T_\si$  is not bounded on $L^2(\rn)$, 
but which satisfies \\
$\sup_x |D_\xi^\al \si(x, \xi)|\leq C_\al |\xi|^{-|\al|}$ and 
$\sup_x \norm{\si(x, \cdot)}{W^{p,n/p}}<\infty$. 
The construction is a minor modification of the 
standard example of a symbol in $S^0_{1, 1}$, which 
is not bounded on $L^2$, see for example \cite{Stein}, page 272. 
We carry out the construction in $n=1$, 
but this can be easily generalized to higher dimensions.

For the given $p>2$,  fix small $0<\de<1/2$, 
so that\footnote{The reason for this choice of $\de$ will become apparent in the proof below.}
 $2+4\de/(1-2\de)<p$. 
Define 
$$
\si(x, \xi):=\suml_{j=8}^\infty \f{e^{-2\pi i 2^j x}}{j^{1/2-\de}}
 \vp(2^{-j} \xi),
$$
where the function $\vp$ is $C^\infty$, -supported 
 in  $1/2\leq |\xi|\leq 3/2$, and $\vp(\eta)=1$ for all $3/4\leq|\eta|\leq 5/4$.  \\
To show the unboundedness of $T_\si$ on $L^2$, 
let us test it against the function 
$$
f_N(x)=\suml_{j=8}^N \f{e^{-2\pi i 2^j x}}{j^{1/2+\de}} f_0(x),
$$
where $f_0$ is a Schwartz function, whose Fourier transform is 
supported in $\{\xi: |\xi|\leq 1/10\}$. Clearly 
$\norm{f_N}{L^2}=(\sum_{j=8}^N j^{1+2\de})^{1/2}\norm{f_0}{L^2}
\leq (\sum_{j=1}^\infty j^{1+2\de})^{1/2}\norm{f_0}{L^2}=C_\de\norm{f_0}{L^2}$, while  
\begin{eqnarray*}
& & T_\si f_N(x)=\sum_{j_1\geq 8, N\geq j_2\geq 8} 
\int \f{e^{-2\pi i 2^{j_1} x}}{j_1^{1/2-\de}} \f{\vp(2^{-j_1} 
\xi)}{j_2^{1/2+\de}} \hat{f_0}(\xi-2^{j_2}) e^{2\pi i \xi x} d\xi.
\end{eqnarray*}
Clearly, by Fourier support considerations the 
terms $j_1\neq j_2$ disappear and we get 
$$
T_\si f_N(x)=(\sum_{N\geq j\geq 8} j^{-1}) f_0(x),
$$
whence $\norm{T_\si}{L^2\to L^2}\gtrsim \ln(N)$, whence $T_\si$ 
{\it is not bounded} on $L^2$. \\ On the other hand, it is clear that for $|\xi|>1$, 
$$
\sup_x |D_\xi^\al \si(x, \xi)|\sim |\xi|^{-|\al|} \ln^{\de-1/2}(|\xi|) 
\leq |\xi|^{-|\al|}.
$$
Finally, to 
estimate $\sup_x \norm{\si(x, \cdot)}{W^{p,1/p}}$, write 
$$
\si=\suml_{s=3}^\infty \suml_{j=2^s}^{2^{s+1}} \f{e^{-2\pi i 2^j x}}{j^{1/2-\de}}
 \vp(2^{-j} \xi)=:\suml_{s=3}^\infty \si^s,
$$
By the convexity of the norms, we have with $\theta: 1/p=\theta/2$, 
$$
\norm{\si^s(x, \cdot)}{W^{p,1/p}}\leq 
\norm{\si^s(x, \cdot)}{H^{1/2}}^\theta
\norm{\si^s(x, \cdot)}{L^\infty}^{(1-\theta)}
$$
It is now easy to compute the norms on the right hand side. We have 
$$
\sup_x \norm{\si^s(x, \cdot)}{H^{1/2}} 
\sim (\sum_{j=2^{s}}^{2^{s+1}} \f{1}{j^{1-2\de}})^{1/2}
\sim 2^{\de s}.
$$
On the other hand, 
$$
\norm{\si^s(x, \cdot)}{L^\infty}\sim  2^{-s(1/2-\de)},
$$
whence $\sup_x \norm{\si^s(x, \cdot)}{W^{p,1/p}}\leq 
2^{s(\de\theta-(1/2-\de)(1-\theta))}$. Clearly, 
such an expression dyadically sums in $s\geq 3$, 
provided $\de\theta<(1/2-\de)(1-\theta)$ or equivalently $p>2+4\de/(1-2\de)$.
\subsection{Proposition \ref{prop:5}: Theorem \ref{theo:4} is sharp}
\begin{proof}(Proposition \ref{prop:5}) 
We construct a sequence of symbols $\si_\de:\rtwo\times \rtwo\to\rone$, 
so that for a fixed Schwartz function $f$
$$
\lim_{\de\to 0+} |T_{\si_\de} f|=|H_* f(x)|=\sup_{u\in \sone} |H_u f(x)|
$$
Since we already know, \cite{LL}, that 
$H_*$ is {\it unbounded} on $L^2(\rtwo)$, 
we should have
\begin{equation}
\label{eq:715}
\limsup_{\de\to 0+} \|T_{\si_\de}\|_{L^2\to L^2}=\infty. 
\end{equation}
In our construction 
$\si_\de$ will depend on $f$, but it is still clear 
that one can achieve \eqref{eq:715}. 
Namely, take a sequence $f_N: \norm{f_N}{L^2(\rtwo)}=1$, 
so that $\|H_* (f_N)\|_{L^2(\rtwo)}\geq N$. 
Then construct $\si_{N, \de}$, so that $\lim_{\de\to 0+} 
|T_{\si_{\de, N}} f_N|=
H_* f_N$. Then clearly, \\
$\limsup_{ N\to \infty, \de\to 0+} 
\|T_{\si_{N,\de}}\|_{L^2\to L^2}=\infty$. 

Now, from the $L^2$ boundedness results of 
Theorem \ref{theo:4} (or rather the lack thereof), we must have 
\begin{equation}
\label{eq:713}
\limsup_{\de\to 0} \suml_l 2^{l/2} 
\sup_x \|P_l^{\xi/|\xi|} \si_\de(x, \cdot)\|_{L^2(\sone)}=\infty.
\end{equation}
On the other hand, we will see that 
$\sup_{x, \xi, \de}|\si_\de(x, \xi)|\leq 1$ and 
\begin{equation}
\label{eq:714}
\sup_{\de, x} \|\si_\de(x, \cdot)\|_{W^{1, 1}(\sone)}<\infty. 
\end{equation}
Note in contrast that (at least heuristically)  \eqref{eq:713} states 
$$
\limsup_{\de\to 0} \sup_x \|\si_\de(x, \cdot)\|_{B^{1/2}_{2, 1}}=\infty
$$
and by the Sobolev embedding estimate on the sphere 
\eqref{eq:bern} (and up  to the usual  Besov spaces 
adjustments at the endpoints), one should have that the quantity 
in \eqref{eq:714} (at least in principle) controls \eqref{eq:713}. Having both 
\eqref{eq:713} and \eqref{eq:714} for a concrete example suggests that 
the conditions imposed in Theorem \ref{theo:4} 
 are extremely tight. 

Let us now describe the construction of $\si_\de$. First of all, 
\begin{eqnarray*}
H_*f(x) = \sup_{u\in\sone} |H_u f(x)| &=& \sup_{u\in\sone}|\int 
sgn(u\cdot \xi/|\xi|) \hat{f}(\xi) e^{2\pi i \xi\cdot x} d\xi|=\\
&=&
|\int 
sgn(u(x)\cdot \xi/|\xi|) \hat{f}(\xi) e^{2\pi i \xi\cdot x} d\xi|,
\end{eqnarray*}
for some measurable function $u(x):\rone\to \sone$. 
Clearly $u(x)$ will depend on the function $f$, see the remarks above after \eqref{eq:715}. \\
Introduce a function $\psi:\psi\in C^\infty, -1\leq \psi(x)\leq 1$, and 
so that $\psi(z)=-1: z\in (-\infty, -1]$, $\psi(z)=1: 
z\in [1, \infty)$. Clearly 
$$
H_* f(x)=\lim_{\de\to 0+} T_{\si_\de} f(x)= \lim_{\de\to 0+}|\int 
\psi\left(\f{u(x)\cdot \xi/|\xi|}{\de}\right) 
\hat{f}(\xi) e^{2\pi i \xi\cdot x} d\xi|, 
$$
that is $\si_\de(x, \xi/|\xi|)=\psi(\f{u(x)\cdot \xi/|\xi|}{\de})$, 
for which we will verify \eqref{eq:714}, while it 
is clearly bounded in absolute value by one. 

We pause 
for a second to comment on the particular form of $T_{\si_\de}$. Note that the function $u(x)$ 
in general will not be smooth\footnote{Note that under some extra smoothness assumptions on $u$, 
Lacey and Li have managed to prove $L^2$ boundedness!} and therefore will not fall under 
the scope of any standard boundedness theory for PDO. Also, note that while the map 
 $\xi\to \si_\de(x, \xi)$ is definitely smooth, its  derivatives are quite large and 
blow up at the important limit $\de\to 0$.  This shows that in order to 
treat maximal operators, build upon singular multipliers (as is the case here), 
one needs the full strength  of Theorems \ref{theo:1}, \ref{theo:4} and beyond.   

Going back to the proof of \eqref{eq:714}, compute 
\begin{eqnarray*}
& & \f{\p \si}{\p_{\xi_1}}= 
\f{\psi'(\f{u(x)\cdot \xi/|\xi|}{\de})}{\de|\xi|^3}\left(u_1(x)\xi_2^2- u_2(x)\xi_1\xi_2\right)= 
\f{\psi'(\f{u(x)\cdot \xi/|\xi|}{\de})\xi_2}{\de|\xi|^2}u(x)\cdot 
(\xi/|\xi|)^{\perp} \\
& & \f{\p \si}{\p_{\xi_2}}=
\f{\psi'(\f{u(x)\cdot \xi/|\xi|}{\de})}{\de|\xi|^3}
\left(u_1(x)\xi_1^2- u_2(x)\xi_1\xi_2\right)= \f{\psi'(\f{u(x)\cdot 
\xi/|\xi|}{\de})\xi_1}{\de|\xi|^2} u(x)\cdot 
(\xi/|\xi|)^{\perp} 
\end{eqnarray*}
Clearly, the supports of both derivatives are in 
$\xi: |u(x)\cdot \xi/|\xi||\leq \de<<1$. Also, on their support, 
$|\nabla \si(x, \xi/|\xi|)|\sim C \de^{-1}$. It follows 
$$
\|\si_\de(x, \cdot)\|_{W^{1, 1}(\sone)}\leq 
\int_{\xi\in \sone: |u(x)\cdot \xi/|\xi|\leq \de}
|\nabla \si_\de(x, \xi)|d\xi\leq C,
$$
where $C$ is independent of $\de$. This was the claim  in \eqref{eq:714}. 
\end{proof}

\section{Appendix}

\subsection{Estimates for Fourier transforms of functions supported on small spherical caps.}
In this section, we present a pointwise estimate for the  kernels of 
multipliers that   restrict the Fourier transform to a small spherical cap. 
\begin{lemma}
\label{le:900}
Let $\theta\in \sn$  and 
$\vp$ is a $C^\infty$ function with $supp \ \vp\subset \{\xi: 1/2\leq |\xi|\leq 2\}$. 
Let also $l>0$ be any integer. Define $K_{l, \theta}$ to be the inverse 
Fourier transform of $\vp(2^l(\xi/|\xi|-\theta))\vp(|\xi|)$, that is 
$$
K_{l, \theta} (x)=\int \vp(2^l(\xi/|\xi|-\theta))\vp(|\xi|) e^{2\pi i x\cdot \xi} d\xi.
$$
Then, for every $N>0$, there exists $C_N$, so that 
\begin{equation}
\label{eq:201}
|K_{l, \theta}(x)|\leq C_N 2^{-l(n-1)}(1+|\dpr{x}{\theta}|)^{-N} 
(1+2^{-l}|x-\dpr{x}{\theta}\theta|)^{-N}.
\end{equation}
That is, in the direction of $\theta$, the function has any polynomial decay, while in the directions transversal to $\theta$, one has decay like $(2^{-l}<x'>)^{-N}$, where $x=\dpr{x}{\theta}\theta+x'$. 
In particular, 
\begin{equation}
\label{eq:202}
\sup_{\theta, l} \int |K_{l, \theta}(x)| dx\leq C_n<\infty,
\end{equation}
where the constant $C_n$ depends on $\|\vp\|_{L^\infty}$ and 
the smoothness properties of $\vp$. 
\end{lemma} 
\begin{proof}
By rotation symmetry, we can assume without loss of generality that $\theta=e_1=(1, 0, \ldots, 0)$. 
Fix  $l$ and drop the subindices for notational convenience. 
We will need to show that for every 
$x=x_1 e_1+x'$, 
\begin{equation}
\label{eq:203}
|K(x)|\leq C_N 2^{-l(n-1)}  <x_1>^{-N} <2^{-l} x'>^{-N}
\end{equation}
First of all,  by support considerations, one has 
$|K(x)|\leq C_n 2^{-l(n-1)}.$ \\
Next, we will show that integration by parts in the variable $\xi_1$ 
yields
\begin{equation}
\label{eq:700} 
K(x)= x_1^{-1} \tilde{K}(x),
\end{equation}
whereas  integration by parts in each of the variables  $\xi_j: j=2, \ldots, n$ yield 
\begin{equation}
\label{eq:701} 
K(x)=(2^{-l} x_j)^{-1} \tilde{K}(x),
\end{equation}
where $\tilde{K}(x)$ is different in each instance, but it has the form 
$$
\tilde{K}(x)=\int \vp_1(2^l(\xi/|\xi|-e_1))\vp_2(|\xi|) e^{2\pi i x\cdot \xi} d\xi.
$$
for some $C^\infty$ 
functions\footnote{As we shall see the functions $\vp_1, \vp_2$ are obtained in a specific way 
from $\vp$ via the operations 
differentiation and  multiplication by monomial.}  
$\vp_1, \vp_2$ with $supp \ \vp_{k}\subset \{\xi: 1/2\leq |\xi|\leq 2\}, k=1,2$. 

That is enough to deduce \eqref{eq:203} and thus Lemma \ref{le:900}. 
Indeed, by iterating  \eqref{eq:700} and  \eqref{eq:701}, one  gets the formula 
$$
K(x)= x_1^{-N_1} (\prod_{j=2}^n (2^{-l} x_j)^{N_j})^{-1} \tilde{K}_{N_1, \ldots, N_n}(x),
$$ 
for any $n$ tuple of integers $(N_1, N_2, \ldots, N_n)$. Combining this representation with the 
estimate $|\tilde{K}_{N_1, \ldots, N_n}(x)|\leq C_{n, N_1, \ldots, N_n} 2^{-l(n-1)}$, one deduces \eqref{eq:203}. 

For \eqref{eq:700}, integration by parts yields 
\begin{eqnarray*}
& & K(x)=-\f{1}{2\pi i x_1} \int \p_{\xi_1} [\vp(2^l(\xi/|\xi|-e_1))\vp(|\xi|)] 
e^{2\pi i x\cdot \xi} d\xi = \\
& & =
-\f{1}{2\pi i x_1}\int  \sum_{j=2}^n 2^l \f{\xi_j^2}{|\xi|^2}
 \p_1\vp(2^l(\xi/|\xi|-e_1)) 
\vp(|\xi|)|\xi|^{-1} e^{2\pi i x\cdot \xi} d\xi+\\
& & 
+\f{1}{2\pi i x_1} \int \sum_{j=2}^n 2^l \f{\xi_j\xi_1}{|\xi|^2} 
\p_j \vp(2^l(\xi/|\xi|-e_1))\vp(|\xi|)|\xi|^{-1} e^{2\pi i x\cdot \xi} d\xi+\\
& &- \f{1}{2\pi i x_1}\int \vp(2^l(\xi/|\xi|-e_1)) \vp'(|\xi|)\f{\xi_1}{|\xi|} 
e^{2\pi i x\cdot \xi} d\xi
\end{eqnarray*}
The third term is clearly in the form $x_1^{-1} \tilde{K}(x)$, by taking into account that $supp \vp\subset \{\xi: 1/2\leq |\xi|\leq 2\}$. \\
The second term above can be rewritten in the form 
\begin{eqnarray*}
& &
\f{1}{2\pi i x_1} \int 
\vp_1(2^l(\xi/|\xi|-e_1))\f{\xi_1}{|\xi|}\vp(|\xi|)|\xi|^{-1} e^{2\pi i x\cdot \xi} d\xi=\\
& & = \f{1}{2\pi i x_1}\int 
[\vp_1(2^l(\xi/|\xi|-e_1))+2^{-l} 
\tilde{\vp}_{1}(2^l(\xi/|\xi|-e_1))]\vp(|\xi|)|\xi|^{-1} e^{2\pi i x\cdot \xi} d\xi=\\
& & = x_1^{-1} \tilde{K}(x),
\end{eqnarray*}
where $\vp_1(\eta)=\sum_{j=2}^n \eta_j \p_{\eta_j}\vp(\eta)$ and 
$\tilde{\vp}_1(\eta)=\eta_1\vp_1(\eta)$. \\
Analogously, one can rewrite the first term of $K(x)$ in the form
$2^{-l} x_1^{-1} \tilde{K}(x)$, i.e. it has an extra decay factor of $2^{-l}$. 
This establishes \eqref{eq:700}. 

For \eqref{eq:701}, we obtain by integration by parts in $\xi_j, 2\leq j\leq n$, 
\begin{eqnarray*}
& & K(x)=\f{1}{2\pi i x_j} \int \p_j\vp(2^l(\xi/|\xi|-e_1))\f{2^l \xi_1 \xi_j}{|\xi|^3} \vp(|\xi|) 
e^{2\pi i x\cdot \xi} d\xi + \\
& &+ \f{1}{2\pi i x_j} \suml_{k\neq j, k=2}^n 
\int \p_k\vp(2^l(\xi/|\xi|-e_1))\f{2^l \xi_k \xi_j}{|\xi|^3} \vp(|\xi|) 
e^{2\pi i x\cdot \xi} d\xi+\\
& &- \f{1}{2\pi i x_j} \int \p_j \vp(2^l(\xi/|\xi|-e_1)) 2^l(\sum_{k\neq j, k=1}^n 
\f{\xi_k^2}{|\xi|^2}) \vp(|\xi|) e^{2\pi i x\cdot \xi} d\xi + \\
& &- \f{1}{2\pi i x_j}\int \vp(2^l(\xi/|\xi|-e_1)) \vp'(|\xi|)\f{\xi_j}{|\xi|} 
e^{2\pi i x\cdot \xi} d\xi
\end{eqnarray*}
By performing similar analysis as in the proof of \eqref{eq:700}, we easily see that the first term above 
is in the form $x_j^{-1} \tilde{K}(x)$, the second and the fourth terms are  in fact even better, since 
they are in the form $2^{-l}x_j^{-1} \tilde{K}(x)$. The third term has two types of terms. Clearly, 
$$
\f{1}{2\pi i x_j} \int \p_j \vp(2^l(\xi/|\xi|-e_1)) 2^l(\sum_{k\neq j, k=2}^n 
\f{\xi_k^2}{|\xi|^2}) \vp(|\xi|) e^{2\pi i x\cdot \xi} d\xi +
$$
is of the form $2^{-l}x_j^{-1} \tilde{K}(x)$, while lastly, 
$$
\f{1}{2\pi i x_j} \int \p_j \vp(2^l(\xi/|\xi|-e_1)) 
2^l \f{\xi_1^2}{|\xi|^2} \vp(|\xi|) e^{2\pi i x\cdot \xi} d\xi
$$
is of the form $2^l x_j^{-1} \tilde{K}(x)$, as is the statement of \eqref{eq:701}. 
\end{proof}
\subsection{$l^p$ functions of cone multipliers}
In this section, we discuss a simple  extension of Lemma \ref{le:900}, 
which is concerned with appropriate 
 $L^p$ bounds 
for $l^p$ functions of such cone multipliers. 
\begin{lemma}
\label{le:sum}
Let $l>>1$ and $\{\theta_m^l\}_{m}$ be a $2^{-l}$ net in $\sn$, so that the family 
$\{\theta\in \sn: 
|\theta_m^l-\theta|\leq 2^{-l}\}_m$ has the finite intersection property. Define 
$$
\widehat{P_m f}(\xi)=\vp_{l, m}(2^l(\xi/|\xi|-\theta_m^l)) \vp(|\xi|) \hat{f}(\xi).
$$
where $\vp_{l, m}$ are as in \eqref{eq:fun}.
Then one has   
\begin{eqnarray}
\label{eq:n1}
& & \|\suml_m P_m g_m\|_{L^p(\rn)}\leq C \left(\suml_m 
\|g_m\|_{L^p}^p\right)^{1/p}\qq \textup{if}\q 1\leq p\leq 2 \\
\label{eq:n2}
& & (\suml_m \|P_m g\|_{L^q(\rn)}^q )^{1/q}\leq  C \|f\|_{L^q}\qq \textup{if}\q 2\leq
q\leq \infty.
\end{eqnarray}
\end{lemma}
\begin{proof} 
Since \eqref{eq:n1} and \eqref{eq:n2} are dual, it will suffice to check \eqref{eq:n2}. Next, the
$L^2$ estimate is trivial by the Plancherel's theorem and the finite intersection property of the
supports of $\vp_{l, m}(2^l(\xi/|\xi|-\theta_m^l))$. Thus, by interpolation it suffices to check 
$$
\supl_m \|P_m g\|_{L^\infty}\leq C\|g\|_{L^\infty}.
$$
But $P_m g(x)=K_{l, \theta_m^l}*g(x)$ and so 
$$
\|P_m g\|_{L^\infty}\leq \|K_{l, \theta_m^l}\|_{L^1}\|g\|_{L^\infty}\leq C \|g\|_{L^\infty}.
$$
where the last inequality follows from \eqref{eq:202}.
\end{proof}

\end{document}